\RequirePackage{fix-cm}

\documentclass[smallextended]{svjour3}       

\smartqed  

\usepackage{graphicx}
\usepackage[margin=1.5in]{geometry}
\usepackage{epstopdf}
\usepackage{color}
\usepackage{amsmath}
\usepackage{multirow}
\usepackage{algorithm}
\usepackage{algorithmic}
\usepackage{subfigure}
\usepackage{booktabs}
\usepackage{dsfont}
\usepackage{amssymb}
\usepackage[latin1]{inputenc}
\usepackage[colorlinks, linkcolor=blue, anchorcolor=blue, citecolor=blue]{hyperref}
\usepackage{caption}
\usepackage{bbding}

\usepackage{pgfplots}
\pgfplotsset{width=10cm}

\begin{document}

\title{Multi-constraint Graph Partitioning Problems Via Recursive Bipartition Algorithm Based on Subspace Minimization Conjugate Gradient Method}

\author{ Wumei Sun \and Hongwei Liu \and  Xiaoyu Wang}

\institute{  Wumei Sun \Envelope \at
              School of Science, Xi'an University of Science and Technology, Xi'an 710054,  China \\
              \email{sunwumei1992@126.com}
         \and Hongwei Liu \at
              School of Mathematics and Statistics, Xidian University, Xi'an 710126,  China\\
              \email{hwliuxidian@163.com}
        \and  Xiaoyu Wang \at
              School of Mathematics and Statistics, Xidian University, Xi'an 710126,  China\\
              \email{wxyhndwl@163.com}
}

\date{Received: date / Accepted: date}
\maketitle

\begin{abstract}
The graph partitioning problem is a well-known NP-hard problem. In this paper, we formulate a 0-1 quadratic integer programming model for the graph partitioning problem with vertex weight constraints and fixed vertex constraints, and propose a recursive bipartition algorithm based on the subspace minimization conjugate gradient method. To alleviate the difficulty of solving the model, the constrained problem is transformed into an unconstrained optimization problem using equilibrium terms, elimination methods, and trigonometric properties, and solved via an accelerated subspace minimization conjugate gradient algorithm. Initial feasible partitions are generated using a hyperplane rounding algorithm, followed by heuristic refinement strategies, including one-neighborhood and two-interchange adjustments, to iteratively improve the results. Numerical experiments on knapsack-constrained graph partitioning and industrial examples demonstrate the effectiveness and feasibility of the proposed algorithm.
\end{abstract}

\section{Introduction}\label{sec1}


The graph partitioning problem is a well-known combinatorial optimization problem, where the objective is to partition the vertices of a graph into two or more subsets under specific constraints while minimizing the sum of edge weights between the subsets. This problem has widespread applications in fields such as parallel computing \cite{Hendrickson2000}, image processing \cite{Shi2000}, very large-scale integrated (VLSI) circuit design \cite{Kahng2011}, and wireless communication \cite{Fairbrother2018}. As Garey et al. demonstrated in \cite{Garey1976}, graph partitioning is an NP-hard problem, making it highly significant to develop more efficient and effective algorithms.


Circuit partitioning is a critical step in the physical design of VLSI circuits \cite{Karypis1998}. It involves dividing a circuit, typically composed of standard cells or other components, into two or more subsets under certain constraints. This process reduces the complexity of VLSI design, meets encapsulation requirements, and enhances performance. Lengauer \cite{Lengauer1990} first proposed modeling circuits using weighted undirected graphs, while Alpert and Kahng \cite{Alpert1995} further formalized the construction of circuit structure diagrams using mathematical models such as weighted undirected and directed graphs. Consequently, solving the circuit partitioning problem can be equated to solving graph or hypergraph partitioning problems. Existing methods for graph partitioning generally fall into two categories: heuristic methods, which aim to quickly compute approximate solutions, and exact methods, which seek optimal partitions.

In the 1960s, Carlson and Nemhauser \cite{Carlson1966} conducted pioneering work on graph partitioning, dividing graphs into up to $k$ subsets without size constraints. Kernighan and Lin \cite{Kernighan1970} introduced the concept of group migration and developed the Kernighan-Lin (KL) algorithm, the first successful heuristic for graph bipartitioning. Subsequent improvements to the KL algorithm emerged, such as the Fiduccia-Mattheyses (FM) algorithm \cite{Fiduccia1982}, which employs a single-point move strategy and bucket sorting to reduce time complexity. While the FM algorithm is efficient, it often converges to local optima and remains limited to single-component optimization. Other heuristic approaches include geometric methods \cite{Gilbert1998}, spectral methods based on the eigenvalues and eigenvectors of the Laplacian matrix \cite{Hendrickson1995a}, multilevel schemes \cite{Hendrickson1995b}, and randomized techniques such as genetic algorithms \cite{Soper2004} and simulated annealing \cite{Kirkpatrick1983}. Popular software implementations of these heuristics include CHACO \cite{Hendrickson1995b}, METIS \cite{Karypis1999}, SCOTCH \cite{Pellegrini1996,Chevalier2008}, Jostle \cite{Walshaw2007}, and KaHiP \cite{Sanders2012,Sanders2013}, among others \cite{Hein2010,Lang2004}. Although these algorithms can process large graphs with millions of vertices efficiently, they do not guarantee optimality and only provide feasible solutions.

Exact methods for graph partitioning include branch-and-cut and branch-and-bound techniques. These methods often rely on complex mechanisms to derive lower bounds, such as linear programming \cite{Brunetta1997}, semidefinite programming \cite{Armbruster2007,Karisch2000}, and quadratic programming \cite{Hager2013}. Brunetta et al. \cite{Brunetta1997} developed a branch-and-cut method based on linear programming relaxation and separation techniques, solving the equipartition problem using edge weights. Karisch et al. \cite{Karisch2000} proposed a branch-and-bound method using semidefinite programming relaxation, which performs well for graph bisection but cannot handle vertex weights or partitions into more than two subsets. Sensen \cite{Sensen2001} introduced a branch-and-bound method based on a generalized lower bound derived from solving a multicommodity flow problem. Hager et al. \cite{Hager2013} applied a branch-and-bound method to continuous quadratic programming, demonstrating that semidefinite programming often yields tighter lower bounds, albeit at high computational cost. Delling et al. \cite{Delling2014} later proposed an exact combinatorial algorithm for minimal graph bisection, employing a branch-and-bound framework with heuristic-based lower bounds. A comprehensive review of graph partitioning applications and recent advances can be found in \cite{Bulu2016}.


Recent research has focused on graph partitioning with vertex weight constraints, particularly the graph partitioning problem under knapsack constraints (GPKC). In GPKC, each vertex is assigned a weight, and subsets must satisfy knapsack constraints. Nguyen \cite{Nguyen2016} combined linear programming relaxation with heuristics to establish upper bounds for GPKC. Semidefinite programming has shown promise in deriving tight lower bounds for quadratic problems with knapsack constraints \cite{Helmberg1996}. Building on this, Wiegele and Zhao \cite{Wiegele2022} successfully applied semidefinite programming relaxation to GPKC, transforming the integer programming model into a semidefinite programming relaxation and solving it using the alternating direction method of multipliers. However, their algorithm required an hour to solve GPKC problems with 500 vertices, making it unsuitable for large-scale problems. Additionally, their approach only considered one-dimensional vertex weights and did not address initial fixed vertex constraints.


In this paper, we propose a recursive bipartitioning approach for graph partitioning with knapsack constraints and initial fixed vertex constraints. Specifically, we formulate the problem as a 0-1 quadratic integer programming model. To address the challenges of solving this nonlinear integer programming problem, we introduce an equilibrium term in the objective function, eliminating vertex weight inequality constraints and relaxing the model to an unconstrained optimization problem. We solve this problem using a subspace minimization conjugate gradient method, followed by random perturbation and an improved hyperplane rounding algorithm to generate an initial feasible solution. Finally, we design a heuristic algorithm to refine the initial partitioning. Numerical experiments demonstrate that our method achieves results comparable to those of \cite{Wiegele2022} in significantly less time.

The paper is organized as follows.  In Section \ref{sec2}, some formulas and models are described, which include how to handle vertex weight constraints and fixed vertex constraints to obtain relaxed unconstrained optimization problems. Then an accelerated subspace minimization conjugate gradient algorithm based on Barzilai-Borwein method(ASMCG\_BB) for unconstrained optimization problems is presented. Section \ref{sec3} gives details of algorithms for graph bipartitioning and recursive bipartitioning methods for multipartitioning of graphs, which include hyperplane rounding algorithms and heuristics. Some numerical experiments are shown in Section \ref{sec4}. Conclusions are given in the last section.


\textbf{Notation} Throughout the paper, $e_{n}$ denotes the vector of all ones of length $n,$ we omit the subscript in case the dimension is clear from the context. The notation $[n]$ denotes the set of integers $\{1,\ldots,n\}.$ The notation $\langle\cdot,\cdot\rangle$ stands for the trace inner product, that is, for any $M, N \in \mathbb{R}^{n\times n},$ define $\langle M, N \rangle := trace(M^{T}N).$ We denote by $diag(M)$ the operation of getting the diagonal entries of matrix $M$ as a vector.

\section{Graph Partitioning Model and Solution of Unconstrained Optimization Problem}\label{sec2}
\subsection{Graph Partitioning Model}
Given an undirected graph $G(V, E)$ with $V = \left\{1,\ldots,n\right\}$ the set of vertices, $E = \left\{\{i, j\} : i, j \in V, i \neq j\right\}$ the set of edges, where $n$ is the number of vertices in graph $G(V, E)$. Each vertex of the graph has a non-negative weight of at least one dimension and the total weight of vertices in each group has the capacity bound. Let $w_{ij}$ be the weight of the edge connecting vertices $i$ and $j$. If there is no edge between vertices $i$ and $j,$ we set $w_{ij}=0.$ For each vertex $i$ and $j,$ we assume that $w_{ii}=0$, $w_{ij}=w_{ji};$ in another word, we are considering an undirected graph without self-loops. The bipartition problem for graphs under the knapsack constraint requires dividing the vertices of the graph into two parts, i.e., set $S$ and set $V\setminus S$, and minimizing the sum of the weights of the edges connecting the vertices in different sets, provided that the sum of the vertex weights of each set does not exceed the upper capacity limit, then the objective of the problem can be written as
\begin{equation}\label{(e2.1)}
 \mathop {\min } \sum\limits_{i\in S, j\in {V\setminus S}} w_{ij}.
\end{equation}
For vertices $i$ and $j,$ if vertice $i\in S,$ let $Y_{i1}=1$ and $ Y_{i2}=0;$ if vertice $j\in {V\setminus S},$ let $Y_{j1}=0$ and $ Y_{j2}=1,$ where 0 and 1 represent the values of two sets. Then,  the multi-constraint graph partitioning problem problem can be formulated as
\begin{align}
 \nonumber \mathop {\min } \ \ &\frac{1}{2}\sum  w_{ij}\left(1-X_{ij}\right)\\\nonumber
 s.t. \ \ &X = YY^{T},\\\label{(e2.2)}
          &Y^{T}B\leq U, \\\nonumber
         &Y_{ik} \in \{0,1\}^{n\times 2},\\\nonumber
         &\forall i \in [n], j \in [n], k \in [2],
\end{align}
where $B$  is the vertex weight matrix and $U$ is the capacity limit matrix. Suppose the vertex has $m$-dimensional weights and the $j$-th dimensional weight of the $i$-th vertex is $B_{ij}$. Then the vertex weight matrix  $B=(B_{ij})$ is a matrix with $n$ rows and $m$ columns. For a bipartitioned graph, the upper capacity limit $U=(U_{ij})$ is a matrix with 2 rows and $m$ columns, where $U_{ij}$ denotes the upper capacity limit of the $j$-th dimensional vertex weight in the $i$-th resource set. In fact, the objective function of the problem \eqref{(e2.2)} is reasonable. For $X = YY^{T},$ if $Y_{i1}\cdot Y_{j1} + Y_{i2}\cdot Y_{j2}=1,$ i.e., $X_{ij}=1,$ then the vertices of edge $(i,j)$ are in the same set and edge $(i,j)\in E$ is not a cut edge; if $Y_{i1}\cdot Y_{j1} + Y_{i2}\cdot Y_{j2}=0,$ i.e., $X_{ij}=0,$ then the vertices of edge $(i,j)$ are not in the same set and edge $(i,j)\in E$ is a cut edge. Naturally, if edge $(i,j)\in E$ is a cut edge, then $(1-X_{ij})=1;$ if edge $(i,j)\in E$ is not a cut edge, then $(1-X_{ij})=0.$ Therefore, the objective function is naturally to find the smallest sum of cut edge weights, i.e., $\mathop {\min } \ \frac{1}{2}\sum\limits_{i,j \in V}  w_{ij}\left(1-X_{ij}\right).$

Let $W=(w_{ij}) \in \mathbb{R}^{n\times n}$ be the weighted adjacency matrix of the graph and introduce the Laplacian matrix $L=diag(We)-W$ associated with the adjacency matrix, where $e$ is an $n$ dimensional column vector with all components 1 and $diag(We) \in \mathbb{R}^{n\times n}$ is a diagonal matrix whose elements on the diagonal are the components of vector $We$. Then, the graph partition problem \eqref{(e2.2)} can be redescribed in the following 0-1 quadratic integer programming form,
\begin{align}
 \nonumber \mathop {\min } \ \ &\frac{1}{2} \left \langle L, YY^{T} \right\rangle\\\label{(e2.3)}
 s.t. \ \ &Y^{T}B\leq U, \\\nonumber
         &Y_{ik} \in \{0,1\}^{n\times 2},\\\nonumber
         &\forall i \in [n],  k \in [2].
\end{align}
Further, let $Y=\left(
                  \begin{array}{cc}
                    y_{1} & y_{2} \\
                  \end{array}
                \right)
\in \mathbb{R}^{n\times 2},$  $y_{1}, y_{2} \in \mathbb{R}^{n\times 1},$ and let $U= \left(\begin{array}{c}
                                                                                                                    U_{1} \\
                                                                                                                    U_{2}
                                                                                                                  \end{array}\right)
\in \mathbb{R}^{2\times m},$  $U_{1}, U_{2} \in \mathbb{R}^{1\times m},$ the problem \eqref{(e2.3)} can be transformed into
\begin{align}
 \nonumber \mathop {\min } \ \ &\frac{1}{2} \left(y_{1}^{T}L y_{1} +  y_{2}^{T}L y_{2}\right)\\\label{(e2.4)}
 s.t. \ \ &y_{1}^{T}B \leq U_1, \\\nonumber
         &y_{2}^{T}B \leq U_2,\\\nonumber
         &y_{1}, y_{2} \in \{0,1\}^{n\times 1}.
\end{align}

In general the solution of nonlinear integer programming is hard. To facilitate the solution of the above quadratic integer programming problem, by replacing condition $y_{1}, y_{2} \in \{0,1\}^{n\times 1}$ with the trigonometric properties, the problem \eqref{(e2.4)} is equivalent to the form of the problem \eqref{(e2.5)} as follows,
\begin{align}
 \nonumber \mathop {\min } \ \ &\frac{1}{2} \left(y_{1}^{T}L y_{1} +  y_{2}^{T}L y_{2}\right)\\\nonumber
 s.t. \ \ &y_{1}^{T}B \leq U_1, \\\label{(e2.5)}
         &y_{2}^{T}B \leq U_2,\\\nonumber
         &y_{1}=\textbf{cos}(\theta), \\\nonumber
         &y_{2}=\textbf{sin}(\theta), \\\nonumber
         & \theta \in \left\{0, \frac{\pi}{2}\right\}^{n\times 1}.
\end{align}
From $y_{1}, y_{2} \in \{0,1\}^{n\times 1}$, it follows that $\textbf{cos}\theta, \textbf{sin}\theta \in \{0,1\}^{n\times 1}$, where
\[
\textbf{cos}(\theta)=\left(\cos(\theta_{1}), \cos(\theta_{2}), \ldots, \cos(\theta_{n})\right)^{T}, \ \textbf{sin}(\theta)=\left(\sin(\theta_{1}), \sin(\theta_{2}), \ldots, \sin(\theta_{n})\right)^{T},
\]
$\theta_{1}, \theta_{2}, \ldots, \theta_{n}$ are all components of the vector $\theta.$
The model \eqref{(e2.5)} is a graph partitioning problem considering only vertex weight constraints, for which several methods have been introduced in Section \ref{sec1}, so we will not repeat them here. In the following, we consider the graph partitioning problem with vertex weight constraints and fixed vertices, i.e., we add the fixed vertex constraints to the model \eqref{(e2.5)}.

For the graph bipartition problem, let the sets of fixed vertices in different sets be $F_{1}$ and $F_{2}$, respectively, and let $F=\left\{F_{1}, F_{2}\right\}$ be the set of all fixed vertices.  Further if vertex $i\in F_{1},$  let $\left(y_{1}\right)_{i}=1$ and $ \left(y_{2}\right)_{i}=0,$ where $\left(y_{1}\right)_{i}$ and $\left(y_{2}\right)_{i}$ denote the $i$-th component of vector $y_{1}$ and vector $y_{2}$ respectively. Similarly, if $i\in F_{2}$, let $\left(y_{1}\right)_{i}=0$ and $ \left(y_{2}\right)_{i}=1,$  then the model \eqref{(e2.5)} with fixed vertices removed by elimination method becomes
\begin{align}
 \nonumber \mathop {\min } \ \ &\frac{1}{2} \left(\bar{y}_{1}^{T}A \bar{y}_{1}+2b_{1}^{T}\bar{y}_{1}+c_{1}\right) + \frac{1}{2} \left(\bar{y}_{2}^{T}A \bar{y}_{2}+2b_{2}^{T}\bar{y}_{2}+c_{2}\right)\\\nonumber
 s.t. \ \ &\bar{y}_{1}^{T}\bar{B} \leq \bar{U}_1, \ \bar{y}_{2}^{T}\bar{B} \leq \bar{U}_2,\\\label{(e2.6)}
         &\bar{y}_{1}=\textbf{cos}(\theta), \ \bar{y}_{2}=\textbf{sin}(\theta),\\\nonumber
         & \theta \in \left\{0, \frac{\pi}{2}\right\}^{\bar{n}\times 1},
\end{align}
where $\bar{y}_{1}$ and $\bar{y}_{2}$ represent the index corresponding to the remaining vertices after removing the fixed vertices;  $A$ is a matrix composed of rows and columns corresponding to non fixed points in the Laplacian matrix $L,$ i.e., $A=\left(L_{ij}\right)$ and $i, j \in V\setminus F;$ $\left(b_{1}\right)_{i}=\sum\limits_{j\in F_{1}} L_{ij},$ $\left(b_{2}\right)_{i}=\sum\limits_{j\in F_{2}} L_{ij},$ $i \in V\setminus F;$ $c_{1}=\sum\limits_{i\in F_{1}}\sum\limits_{j\in F_{1}}L_{ij},$  $c_{2}=\sum\limits_{i\in F_{2}}\sum\limits_{j\in F_{2}}L_{ij};$
$\bar{B}$ is the vertex weight matrix after eliminating the row of fixed vertices; $\bar{U}_1$ and $\bar{U}_2$  are the remaining resource upper limits after the original resource upper limit $ U_{1}$ and  $ U_{2}$ subtracts the vertex weights of fixed tops in sets $F_{1}$ and $F_{2}$ respectively; and $\bar{n}$ is the number of remaining vertices after removing fixed vertices.

It is easy to see that variables in model \eqref{(e2.6)}  are discrete variables taking 0 and $\frac{\pi}{2},$ which means it is relatively difficult to solve and cannot be solved by various excellent nonlinear programming methods.  Motivated by this, we want to do some processing on the model to facilitate the solution, so we relax the value range of the variables in the model \eqref{(e2.6)} to obtain a continuous optimization model as
\begin{align}
 \nonumber \mathop {\min } \ \ &\frac{1}{2} \left(\bar{y}_{1}^{T}A \bar{y}_{1}+2b_{1}^{T}\bar{y}_{1}+c_{1}\right) + \frac{1}{2} \left(\bar{y}_{2}^{T}A \bar{y}_{2}+2b_{2}^{T}\bar{y}_{2}+c_{2}\right)\\\nonumber
 s.t. \ \ &\bar{y}_{1}^{T}\bar{B} \leq \bar{U}_1, \ \bar{y}_{2}^{T}\bar{B} \leq \bar{U}_2,\\\label{(e2.7)}
         &\bar{y}_{1}=\textbf{cos}(\theta), \ \bar{y}_{2}=\textbf{sin}(\theta),\\\nonumber
         & \theta \in \left\{\left[0, \frac{\pi}{2}\right]\right\}^{\bar{n}\times 1}.
\end{align}
Each component of vector $\theta$ in model \eqref{(e2.7)} is taken to be on a continuous interval $\left[0, \frac{\pi}{2}\right].$ The constrained optimization problem of model \eqref{(e2.7)} is usually solved directly by the projected gradient method, the proximal gradient method and the ADMM method, etc. Another type of method is to do some processing of the constraints to transform the problem \eqref{(e2.7)} into an unconstrained optimization problem for solution, such as further relaxing the components of $\theta$ to the interval $\left[0, 2\pi\right]$, while using penalty function method, Lagrangian method or generalized Lagrangian method to penalize the vertex weight inequality constraints to the objective function, and then transform it into an unconstrained problem, which can be finally solved by using various nonlinear programming algorithms for solving unconstrained optimization problems.

In this paper, the model \eqref{(e2.6)}  is also treated by transforming it into an unconstrained problem to be solved, but the treatment of the vertex weight inequality constraints differs from that of the penalty function method. Due to the existence of vertex weight inequality constraints, the nonlinear programming model \eqref{(e2.6)} becomes complex, so we consider removing the vertex weight inequality constraints. However, in practical applications, if the vertex weight inequality constraint is directly removed, some unpredictable problems will occur, because $w_{ij}$ is often considered as a non negative number. If the problem \eqref{(e2.6)} is minimized without fixed vertex constraints, it may occur that the components of the optimal solution vector $\theta$ obtained are all 0 or all $\frac{\pi}{2}$, then the minimum objective value obtained is 0, that is, the vertices are all divided in the same set, further resulting in the weights of these vertices and exceeding the capacity limit. To say the least, even if the fixed point constraint is considered, a similar situation may occur when there are few fixed points. At this time, the ideal division will not be obtained, which is contrary to our solution idea. Therefore, we hope to find a method that can ensure the quality of the optimal solution while removing the complex inequality constraints. To avoid the above situation and reduce the resource waste caused by partitioning, we consider removing the vertex weight inequality constraint and adding an equilibrium term to the objective function of the problem \eqref{(e2.5)} to control the equilibrium of the partitioning results. The reason why the equilibrium term is added to the model \eqref{(e2.5)} is that there are also fixed vertices in the equilibrium term. Then, the problem \eqref{(e2.5)} becomes
\begin{align}
 \nonumber \mathop {\min } \ \ &\frac{1}{2} \left(y_{1}^{T}L y_{1} +  y_{2}^{T}L y_{2}\right) + \frac{1}{2} \rho h(y_{1},y_{2})\\\nonumber
 s.t. \ \  &y_{1}=\textbf{cos}(\theta), \\\label{(e2.8)}
         &y_{2}=\textbf{sin}(\theta), \\\nonumber
         & \theta \in \left\{0, \frac{\pi}{2}\right\}^{n\times 1},
\end{align}
where $\rho h(y_{1},y_{2})$ is the equilibrium term, $\rho$ is the equilibrium factor. The selection of the equilibrium term is crucial to the solution of the problem. A suitable equilibrium term needs to satisfy two requirements, one of which is that it is a non negative term as small as possible, the other is that it should ensure that the partition result is as balanced as possible, for example, the equilibrium term can be $\rho \left(e^{T}y_{1}-e^{T}y_{2}\right)^{2}$ or $\rho \left(\left(e^{T}y_{1}\right)^2 + \left(e^{T}y_{2}\right)^2\right).$  It is easy to find that for the first equilibrium term, when it is as small as possible, $\rho$ is a constant, meaning that the number of nonzero elements in $y_{1}$ and $y_{2}$ is approximately equal, that is, the number of vertices divided into two sets is relatively close, which at this point also means that the partition  is relatively balanced and therefore it is reasonable. For the second equilibrium term, to make $\left(e^{T}y_{1}\right)^2 + \left(e^{T}y_{2}\right)^2$ as small as possible, it can only be satisfied when $e^{T}y_{1}$ and $e^{T}y_{2}$ are approximately equal, thus also giving  a balanced partition result. In particular, $\left(e^{T}y_{1}\right)^2 + \left(e^{T}y_{2}\right)^2$ reaches a minimum when $e^{T}y_{1} = e^{T}y_{2},$ because $e^{T}y_{1} + e^{T}y_{2}$ is the number of all vertices.
After several experiments, it is found that the numerical performance of choosing the equilibrium term $\rho \left(\left(e^{T}y_{1}\right)^2 + \left(e^{T}y_{2}\right)^2\right)$ will be better, so in this paper, we choose the equilibrium term $\rho h(y_{1},y_{2})=\rho \left(\left(e^{T}y_{1}\right)^2 + \left(e^{T}y_{2}\right)^2\right)$. Then, we get the model of adding equilibrium term as
\begin{align}
 \nonumber \mathop {\min } \ \ &\frac{1}{2} \left(y_{1}^{T}L y_{1} +  y_{2}^{T}L y_{2}\right) + \frac{1}{2} \rho \left(\left(e^{T}y_{1}\right)^2 + \left(e^{T}y_{2}\right)^2\right)\\\nonumber
 s.t. \ \  &y_{1}=\textbf{cos}(\theta), \\\label{(e2.9)}
         &y_{2}=\textbf{sin}(\theta), \\\nonumber
         & \theta \in \left\{0, \frac{\pi}{2}\right\}^{n\times 1}.
\end{align}
Similar to the processing of \eqref{(e2.5)} - \eqref{(e2.6)}, after eliminating the fixed vertices using the elimination method, the model \eqref{(e2.9)}  becomes
\begin{align}
 \nonumber \mathop {\min } \ \ &\frac{1}{2} \left(\bar{y}_{1}^{T}\bar{A} \bar{y}_{1}+2\bar{b}_{1}^{T}\bar{y}_{1}+\bar{c}_{1}\right) + \frac{1}{2} \left(\bar{y}_{2}^{T}\bar{A} \bar{y}_{2}+2\bar{b}_{2}^{T}\bar{y}_{2}+\bar{c}_{2}\right)\\\label{(e2.10)}
 s.t. \ \ &\bar{y}_{1}=\textbf{cos}(\theta), \ \bar{y}_{2}=\textbf{sin}(\theta),\\\nonumber
         & \theta \in \left\{0, \frac{\pi}{2}\right\}^{\bar{n}\times 1},
\end{align}
where $\bar{A}=A+\rho \cdot ee^{T},$ $\left(\bar{b}_{1}\right)_{i}=\left(b_{1}\right)_{i}+\rho \left(\sum\limits_{i\in F_{1}} e_{i}\right)e,$ $\left(\bar{b}_{2}\right)_{i}=\left(b_{2}\right)_{i}+\rho \left(\sum\limits_{i\in F_{2}} e_{i}\right)e,$ $\bar{c}_{1}=c_{1}+\rho \left(\sum\limits_{i\in F_{1}} e_{i}\right)^{2},$ $\bar{c}_{2}=c_{2}+\rho \left(\sum\limits_{i\in F_{2}} e_{i}\right)^{2},$ and other quantities are the same as in \eqref{(e2.6)}.
Similar to \eqref{(e2.7)}, we relax the values of the components of $\theta$ to the interval $\left[0,\frac{\pi}{2}\right],$ i.e., $\theta \in \left\{\left[0, \frac{\pi}{2}\right]\right\}^{\bar{n}\times 1},$ and subsequently further relax the values of the components of $\theta$ to the interval $\left[0,2\pi\right]$ in order to remove the constraint about $\theta$. Since the cosine and sine functions are periodic functions with a period of $2\pi$, it is equivalent to restrict the values of the variables to the interval $\left[0,2\pi\right]$ and to remove the constraints on the variables. At this point, the unconstrained optimization problem obtained after the relaxation of the problem \eqref{(e2.10)} is
\begin{align}\label{(e2.11)}
\mathop {\min }\limits_{\theta\in \mathbb{R}^{\bar{n}}} \ \ \frac{1}{2} \left({\textbf{cos}(\theta)}^{T}\bar{A} {\textbf{cos}(\theta)}+2\bar{b}_{1}^{T}{\textbf{cos}(\theta)}+\bar{c}_{1}\right) + \frac{1}{2} \left({\textbf{sin}(\theta)}^{T}\bar{A} {\textbf{sin}(\theta)}+2\bar{b}_{2}^{T}{\textbf{sin}(\theta)}+\bar{c}_{2}\right).
\end{align}
Obviously, the problem \eqref{(e2.11)} is an easy-to-solve unconstrained optimization problem, and there are many nonlinear programming methods to solve it, and the conjugate gradient method is one of the most effective methods. In this paper, we use the accelerated subspace minimization conjugate gradient method to solve the more relaxed unconstrained optimization model \eqref{(e2.11)} with the vertex weight constraint removed, and then apply random perturbation to the obtained solution, and add the feasibility determination of the obtained solution in the process of random perturbation to ensure that the obtained solution satisfies the vertex weight constraint, and then further adjust the vertices to finally obtain the partition result. In the next subsection, we will introduce an accelerated subspace minimization conjugate gradient algorithm based on BB method for solving unconstrained optimization problems.

\subsection{ASMCG Algorithm for Unconstrained Optimization Problems}
The subspace minimization conjugate gradient method is a very interesting method, which shows particularly effective numerical performance in solving unconstrained optimization problems \cite{Andrei2014,Sun21,Sun22,Liu19,Zhao21}. For the subspace minimization conjugate gradient method based on the BB method, Sun and Liu \cite{Sun21} proposed some useful acceleration conditions and then gave an accelerated subspace minimization conjugate gradient method(ASMCG\_BB), in which the acceleration parameter obtained by the quadratic interpolation function was used to improve the step size. Moreover, the convergence of the algorithm is guaranteed under mild assumptions, and the numerical results also show that ASMCG\_BB is effective and feasible. In view of this, we consider solving the above unconstrained optimization problem \eqref{(e2.11)} by using ASMCG\_BB method, which has the following iterative form:
\begin{align}\label{(e2.12)}
\theta_{k+1} =\theta_{k} +\alpha_{k}d_{k},
\end{align}
where $\alpha_{k}$ is step size obtained by a nonmonotone generalized Wolfe line search, and $d_{k}$ is the search direction in the form of:\\
\begin{align}\label{(e2.13)}
d_{k} = u_{k} g_{k} +v_{k} s_{k-1},
\end{align}
where $s_{k-1}=\theta_{k}-\theta_{k-1},$ $g_{k}$ is the gradient of the objective function at $\theta_{k},$ $u_{k}$ and $v_{k}$ are parameters. In other words, the search direction  $d_{k}$ is computed by minimizing the  following problem
\begin{align}\label{(e2.14)}
\min\limits_{d_{k} \in \Omega_{k}} \ \   g_{k }^Td_{k} + \frac{1}{2}{d_{k}^T}{B_{k}}d_{k},
\end{align}
in a two-dimensional subspace $\Omega_{k}=Span\{g_{k},s_{k-1}\},$ where ${{B}_{k}} \in \mathbb{R}^{n \times n}$ is an approximation of Hessian matrix, which is positive definite and symmetric. Substituting \eqref{(e2.13)} to \eqref{(e2.14)} and utilizing the standard secant equation $ B_{k}s_{k-1}=y_{k-1}, $ where $y_{k-1}=g_{k}-g_{k-1},$ then we have:
\begin{align}\label{(e2.15)}
\min\limits_{u,v \in \mathbb{R}} {\left( \begin{array}{c}
                                     \| g_{k} \|^{2} \\
                                     g_{k}^{T}s_{k-1}
                                   \end{array}
\right )} ^{T} {\left( \begin{array}{c}
                         u_{k} \\
                         v_{k}
                       \end{array}\right )} + \frac{1}{2}{\left( \begin{array}{c}
                         u_{k} \\
                         v_{k}
                       \end{array}\right )}^{T} {\left(\begin{array}{cc}
                                                         \rho_{k} &  g_{k}^{T}y_{k-1} \\
                                                         g_{k}^{T}y_{k-1} & s_{k-1}y_{k-1}
                                                       \end{array}\right )}{\left( \begin{array}{c}
                         u_{k} \\
                         v_{k}
                       \end{array}\right ).}
\end{align}
As an estimate of $g^{T}_{k}B_{k}g_{k}$, $\rho_{k}$ is a very significant parameter, whose different values will lead to different search directions. Different $u_{k}$ and $v_{k}$ are obtained according to different quadratic approximation models, and the search directions are further obtained in the following different cases.

Case I

If condition (2.11) in \cite{Sun21} hold, set $u_{k}=0$ in \eqref{(e2.13)} and substitute \eqref{(e2.13)} into \eqref{(e2.15)}, then we have
\begin{align}\label{(e2.16)}
\left( \begin{array}{c}
u_k^1\\
v_k^1
\end{array} \right) = \left( \begin{array}{c}
0\\
- \frac{{g_k^{ T}{s_{k - 1}}}}{{s_{k - 1}^{ T}{y_{k - 1}}}}
\end{array} \right).
\end{align}

Case II

If condition (2.14) in \cite{Sun21} hold, set $\rho_{k}=\frac{3}{2} \frac{{{{\left\| {{y_{k-1}}} \right\|}^2}}}
{{s_{k-1}^T{y_{k-1}}}}{\left\| {{g_{k}}} \right\|^2}.$ Then by minimizing \eqref{(e2.15)}, we get
\begin{align}\label{(e2.17)}
\left( \begin{array}{c}
u_k^2\\
v_k^2
\end{array} \right) = \left( \begin{array}{c}
\frac{1}{\Delta_{k}} \left( g^{T}_{k}y_{k-1} g^{T}_{k}s_{k-1} - s^{T}_{k-1}y_{k-1}\|g_{k}\|^{2} \right)\\
\frac{1}{\Delta_{k}} \left( g^{T}_{k}y_{k-1}\|g_{k}\|^{2}-\rho_{k}g^{T}_{k}s_{k-1} \right)
\end{array} \right),
\end{align}
where
\begin{equation}
\Delta_{k} = {\left|\begin{array}{cc}
                                                         \rho_{k} &  g_{k}^{T}y_{k-1} \\
                                                         g_{k}^{T}y_{k-1} & s_{k-1}y_{k-1}
                                                       \end{array}\right |}  >0.
\end{equation}

Case II

If condition (2.19) in \cite{Sun21} hold, set
\begin{align*}
\rho_{k}=\| g_{k} \|^{2} + \frac{(g^{T}_{k}y_{k-1})^{2}}{s^{T}_{k-1}y_{k-1}}.
\end{align*}
Substituting $\rho_{k}$ into \eqref{(e2.15)}, we obtain
\begin{align}\label{(e2.19)}
\left( \begin{array}{c}
u_k^3\\
v_k^3
\end{array} \right) = \left( \begin{array}{c}
- 1 + \frac{{g_k^{T}{y_{k - 1}}g_k^{T}{s_{k - 1}}}}{{s_{k - 1}^{T}{y_{k - 1}}\parallel {g_k}{\parallel ^2}}}\\
(1 - \frac{{g_k^{T}{y_{k - 1}}g_k^{T}{s_{k - 1}}}}{{s_{k - 1}^{T}{y_{k - 1}}\parallel {g_k}{\parallel ^2}}})\frac{{g_k^{T}{y_{k - 1}}}}{{s_{k - 1}^{T}{y_{k - 1}}}} - \frac{{g_k^{T}{s_{k - 1}}}}{{s_{k - 1}^{T}{y_{k - 1}}}}
\end{array} \right).
\end{align}

Case IV

If none of the previous three cases holds, we choose the negative gradient direction, i.e.
\begin{align*}
{d_k} = -{g_k}.
\end{align*}

In summary, the search direction of ASMCG\_BB method can be described as
\begin{align}\label{(e2.20)}
{d_k} = \left\{ \begin{array}{l}
u_k^1{g_k} + v_k^1{s_{k - 1}},\;\;\;\; \text{Case I},\\
u_k^2{g_k} + v_k^2{s_{k - 1}},\;\;\;\; \text{Case II},\\
u_k^3{g_k} + v_k^3{s_{k - 1}},\;\;\;\; \text{Case III},\\
\;\;\;\;\;\;\; - {g_k}, \;\;\;\;\;\;\;\;\;\;\;\;\;\;\text{Case IV}.
\end{array} \right.
\end{align}

For the completeness of the paper, we briefly introduce the algorithmic framework of ASMCG\_BB method below. Some details and  conditions have been written in great detail in the literature \cite{Sun21}, so we will not repeat them here.

\vspace{0.6cm}
\begin{algorithm}
\caption{ASMCG\_BB}\label{algorithm.1}
\noindent  \textbf{Step 0.} Initialization, give the initial point $\theta_{0}$ of the algorithm and determine the values of various \\
\hspace*{1.2cm} parameters, set $k:=0.$\\
\textbf{Step 1.} If the termination condition $\|g_{k}\|_{\infty} \leq \varepsilon$  holds, stop iteration; otherwise, go to \textbf{Step 2.} \\
\textbf{Step 2.} Compute the search direction $d_{k}$ by \eqref{(e2.20)}.\\
\textbf{Step 3.} Determine the corresponding initial step size $\alpha^{0}_{k}$ according to the different iteration direction \\
\hspace*{1.2cm} $d_{k}$ in \textbf{Step 2.}\\
\textbf{Step 4.} Determine a step size  $\alpha_{k}$ satisfying the generalized nonmonotone Wolfe line search with initial \\ \hspace*{1.2cm} step size $\alpha^{0}_{k}.$\\
\textbf{Step 5.} Compute the trial iteration, and judge whether the gradient at the trial iteration satisfies the\\
\hspace*{1.2cm} termination condition $\|g_{k}\|_{\infty} \leq \varepsilon,$ If yes, stop iteration;  otherwise, go to \textbf{Step 6.}\\
\textbf{Step 6.} Acceleration procedure:\\
\hspace*{1.2cm} \textbf{(1)} Judge whether the acceleration criteria are satisfied. If yes, update the acceleration \\
\hspace*{1.8cm} iteration formula with $\theta_{k+1} =\theta_{k} +\eta_{k}\alpha_{k}d_{k},$ where $\eta_{k}$ is the acceleration parameter   \\
\hspace*{1.8cm} obtained by quadratic interpolation function; \\
\hspace*{1.2cm} \textbf{(2)} Calculate the corresponding function value and gradient value of iteration point $\theta_{k+1},$  \\
\hspace*{1.8cm} then go to \textbf{Step 8.}; otherwise, go to \textbf{Step 7.}\\
\textbf{Step 7.} Update the iteration formula with \eqref{(e2.12)}, calculate the corresponding function value and  \\
\hspace*{1.2cm}gradient value of iteration point $\theta_{k+1}.$ \\
\textbf{Step 8.} Update restart conditions and parameters in line search. Set $k:=k+1$ and go to \textbf{Step 1.} \\
\end{algorithm}
\vspace{0.3cm}

\section{Bipartition and Multipartition Algorithms for Graphs}\label{sec3}
In Section \ref{sec2}, we gave the graph partitioning model with vertex weights and fixed vertex constraints, and described in detail how to relax the graph partitioning model to an unconstrained optimization problem that is easy to solve, and finally gave an accelerated subspace minimization conjugate gradient method to solve it. In this section, we will randomly perturb the solution obtained in Section \ref{sec2} and add the feasibility determination of the solution with respect to the vertex weight constraint, and then further adjust the vertices using one-neighborhood, pairwise tuning and other techniques to finally give satisfactory two-partition or multi-partition results.

\subsection{Heuristic algorithm for bipartition of graphs}
By using Algorithm 1, we obtain an optimal solution vector $\theta^{*} \in R^{\bar{n}}$ of the unconstrained optimization problem \eqref{(e2.11)} after relaxation, and accordingly matrix $ \left(
                                                             \begin{array}{cc}
                                                               \textbf{cos} (\theta^{*}) & \textbf{sin} (\theta^{*}) \\
                                                             \end{array}
                                                           \right) ^{T} \in R^{2\times \bar{n}}$
can be obtained. It is obvious that the components in vectors $\textbf{cos} (\theta^{*})$ and $\textbf{sin} (\theta^{*})$ take values in the continuous interval $\left[-1, 1\right],$ however, the components of the solution vectors $\bar{y}_{1}$ and $\bar{y}_{2}$ representing the vertices of the graph take values in the discrete points 0 and 1. In order to divide the vertices of the graph, we need to correspond the $\textbf{cos} (\theta^{*})$ and $\textbf{sin} (\theta^{*})$  in \eqref{(e2.11)} to the discrete $\bar{y}_{1}$ and $\bar{y}_{2}$ whose components are integers in \eqref{(e2.10)}.
Then a hyperplane is randomly selected and used to divide the solution vector representing the vertices of the undirected graph into two parts to obtain a preliminary partition result for the graph vertices. Then, we add weight constraints to these vertices of the initial partition, and judge whether the vertex weight sum exceeds the capacity limit, if it is within the resource limit, then this partition result is feasible, otherwise we need to further adjust the partition result.

From an algorithmic point of view, this method of obtaining the preliminary partitioning results of the graph is actually equivalent to the hyperplane rounding method. First, a set of samples $\left( \gamma_{1}, \gamma_{2}, \ldots, \gamma_{p}, \right)$ that obey uniform distribution on $0-2\pi$ is selected, and the vector $ \left(
                                                                           \begin{array}{cc}
                                                                             \cos(\gamma_{j}) & \sin(\gamma_{j}) \\
                                                                           \end{array}
                                                                         \right)
^{T}, j\in [p],$ consisting of the cosine and sine functions corresponding to each sample is selected to be the normal vector $\beta_{j}$ of the hyperplane, i.e., the hyperplane is $H=\left\{v | \beta_{j}^{T} v =0, j\in [p] \right\}$. Then, the matrix composed of normal vectors of $p$ different hyperplanes is
\begin{align}\label{(e3.1)}
\beta = \left( \beta_{1}, \beta_{2}, \ldots, \beta_{p}\right) = \left(
                                                   \begin{array}{cccc}
                                                     \cos(\gamma_{1}) & \cos(\gamma_{2}) & \ldots & \cos(\gamma_{p}) \\
                                                     \sin(\gamma_{1}) & \sin(\gamma_{2}) & \ldots & \sin(\gamma_{p}) \\
                                                   \end{array}
                                                 \right) \in R^{2\times p}.
\end{align}
Let the matrix $s$ consisting of cosine function and sine function corresponding to $\theta^{*}$ in \eqref{(e2.11)} be
\begin{align}\label{(e3.2)}
s = \left( s_{1}, s_{2}, \ldots, s_{\bar{n}}\right) = \left(
      \begin{array}{cccc}
        \cos (\theta_{1}^{*}) & \cos (\theta_{2}^{*}) & \ldots & \cos (\theta_{\bar{n}}^{*}) \\
        \sin (\theta_{1}^{*}) & \sin (\theta_{2}^{*}) & \ldots & \sin (\theta_{\bar{n}}^{*}) \\
      \end{array}
    \right) \in R^{2\times \bar{n}},
\end{align}

Then the value of approximate solution vectors $\bar{y}_{1}$ and $ \bar{y}_{2} $ in \eqref{(e2.10)} is rounded according to the value of $ s^{T}\beta_{j} $, that is, if $s_{i}^{T}\beta_{j} \geq 0 (i\in [\bar{n}], j\in [p]),$ take the $i$-th component of vector $\bar{y}_{1}$ as $\left(\bar{y}_{1}\right)_{i}=1$ and the $i$-th component of vector $\bar{y}_{2}$ as $\left(\bar{y}_{2}\right)_{i}=0;$ if $s_{i}^{T}\beta_{j} < 0,$ take the $i$-th component of vector $\bar{y}_{1}$ as $\left(\bar{y}_{1}\right)_{i}=0$ and the $i$-th component of vector $\bar{y}_{2}$ as $\left(\bar{y}_{2}\right)_{i}=1$. Here $\left(\bar{y}_{1}\right)_{i},\left(\bar{y}_{2}\right)_{i} \in \left\{0, 1\right\}$ represent the vertice $i$, $0$ and $1$ represent the values of two sets, namely, $S=\left\{ i: \left(\bar{y}_{1}\right)_{i}=1, \left(\bar{y}_{2}\right)_{i}=0 \right\}, V\backslash S=\left\{ i: \left(\bar{y}_{1}\right)_{i}=0, \left(\bar{y}_{2}\right)_{i}=1 \right\}$. In this way, we obtain a preliminary partitioning result of the graph without vertex weight constraint. The feasibility determination is then performed based on this preliminary partitioning result, which leads to a feasible solution that satisfies both the problem \eqref{(e2.10)} and the following vertex weight inequality constraint condition\\
\begin{align}\label{(e3.3)}
\bar{y}_{1}^{T}\bar{B} \leq \bar{U}_1, \ \bar{y}_{2}^{T}\bar{B} \leq \bar{U}_2.
\end{align}

Then the approximate optimal solution of the graph partitioning problem is obtained by using heuristic strategies such as one-neighborhood adjustment and vertex swapping adjustment. In view of this, the hyperplane rounding algorithm for the graph bipartition problem is given.

\vspace{0.5cm}
\begin{algorithm}
	\caption{The hyperplane rounding algorithm for the bipartition problem of graphs}
\noindent
\textbf{Data:} Vertex weight matrix $\bar{B},$  resource upper limits $\bar{U}_1$ and $\bar{U}_2;$ \\
\textbf{Input:} The optimal solution $\theta^{*}$ of unconstrained optimization problem \eqref{(e2.11)};\\
\textbf{Output:} Vectors $\bar{y}_{1}^{*}$ and $\bar{y}_{2}^{*}$ representing vertex partitions in different sets, partitions $P_{1}$ and $P_{2}$.\\  \indent
\textbf{Step 0.} Select the uniformly distributed $p$-dimensional vector $\gamma$ on the unit circle, and calculate the \\  \indent \hspace*{1.2cm} matrix $\beta$ composed of normal vectors of different hyperplanes according to \eqref{(e3.1)}.\\  \indent
\textbf{Step 1.} Generate the $p$ group approximate solutions $\bar{y}_{1}$ and $\bar{y}_{2}$ of problem \eqref{(e2.10)}.\\  \indent
\textbf{Step 2.} Calculate the objective function $f\left(\bar{y}_{1},\bar{y}_{2}\right)$ of problem \eqref{(e2.10)} corresponding to the $p$-group\\  \indent
\hspace*{1.2cm} approximate solution in \textbf{Step 1.}; find the approximate solutions corresponding to the smallest \\  \indent \hspace*{1.2cm} $\bar{p} (\bar{p} \leq p)$ sets of function values $f_{min}\left(\bar{y}_{1},\bar{y}_{2}\right)$ among  $p$ sets of function values noted as ${\bar{y}_{1}}{_{min}}$ \\  \indent
\hspace*{1.2cm} and ${\bar{y}_{2}}{_{min}}.$\\  \indent
\textbf{Step 3.} Judge whether the constraint conditions of vertex weight inequality are satisfied:\\  \indent
\hspace*{1.2cm} If ~~ ${\bar{y}_{1}}{_{min}}^{T}\bar{B} \leq \bar{U}_1$ and ${\bar{y}_{2}}{_{min}}^{T}\bar{B} \leq \bar{U}_2$  \\ \indent
\hspace*{1.8cm} Output $\bar{y}_{1}^{*} :={\bar{y}_{1}}{_{min}}$ and $\bar{y}_{2}^{*} :={\bar{y}_{2}}{_{min}};$ go to \textbf{Step 4.}\\ \indent
\hspace*{1.2cm} else \\ \indent
\hspace*{1.8cm} If ~~ $\bar{y}_{1}^{T}\bar{B} \leq \bar{U}_1$ and $\bar{y}_{2}^{T}\bar{B} \leq \bar{U}_2$  \\ \indent
\hspace*{2.5cm} \textbf{(1)} Calculate the objective function $f\left(\bar{y}_{1},\bar{y}_{2}\right)$ of problem \eqref{(e2.10)}; \\ \indent
\hspace*{2.5cm} \textbf{(2)} Find the approximate solutions corresponding to the smallest function values \\  \indent \hspace*{2.5cm} $f_{min}\left(\bar{y}_{1},\bar{y}_{2}\right)$  noted as ${\bar{y}_{1}}{_{min}}$ and ${\bar{y}_{2}}{_{min}};$\\  \indent
\hspace*{2.5cm} \textbf{(3)} Output $\bar{y}_{1}^{*} :={\bar{y}_{1}}{_{min}},$ and $\bar{y}_{2}^{*} :={\bar{y}_{2}}{_{min}};$ go to \textbf{Step 4.} \\ \indent
\hspace*{1.8cm} else \\ \indent
\hspace*{2.5cm} Adjust $\bar{y}_{1}$ and $\bar{y}_{2}$ with one-neighborhood and two-interchanges  to get $\bar{\bar{y}}_{1}$ and $\bar{\bar{y}}_{2};$ \\ \indent
\hspace*{2.5cm} If ~~ $\bar{\bar{y}}_{1}^{T}\bar{B} \leq \bar{U}_1$ and $\bar{\bar{y}}_{2}^{T}\bar{B} \leq \bar{U}_2$  \\ \indent
\hspace*{3.2cm} \textbf{(1)} Calculate the objective function $f\left(\bar{\bar{y}}_{1},\bar{\bar{y}}_{2}\right)$ of problem \eqref{(e2.10)}; \\ \indent
\hspace*{3.2cm} \textbf{(2)} Find the approximate solutions corresponding to the smallest function \\ \indent
\hspace*{3.2cm} values  $f_{min}\left(\bar{\bar{y}}_{1},\bar{\bar{y}}_{2}\right)$  noted as ${\bar{\bar{y}}_{1}}{_{min}}$ and ${\bar{\bar{y}}_{2}}{_{min}};$ \\  \indent
\hspace*{3.2cm} \textbf{(3)} Output $\bar{y}_{1}^{*} :={\bar{\bar{y}}_{1}}{_{min}},$ and $\bar{y}_{2}^{*} :={\bar{\bar{y}}_{2}}{_{min}};$ go to \textbf{Step 4.} \\ \indent
\hspace*{2.5cm} else \\ \indent
\hspace*{3.2cm} Add resources and update resource upper limits $\bar{U}_1$ and $\bar{U}_2;$ go to \textbf{Step 3.}\\ \indent
\hspace*{2.5cm} end \\ \indent
\hspace*{1.8cm} end \\ \indent
\hspace*{1.2cm} end \\ \indent
\textbf{Step 4.} Add indicators corresponding to fixed vertices and combine $\bar{y}_{1}^{*}$ and $\bar{y}_{2}^{*}$ to output partitions $P_{1}$ \\ \indent
\hspace*{1.2cm} and $P_{2}.$
\end{algorithm}
\vspace{0.5cm}

\emph{\textbf{Remark 1.}} In algorithm 2, if $\forall i\in [n]$ represents all vertices to be divided, then partition $P_{1}=\left\{ i| \left(y_{1}\right)_{i}=1, \left(y_{2}\right)_{i}=0 \right\}$ and $P_{2}=\left\{ i| \left(y_{1}\right)_{i}=0, \left(y_{2}\right)_{i}=1 \right\}.$   Two vertex adjustment strategies, one-neighborhood and two-interchanges, are a heuristic algorithmic idea that can be seen as greedy algorithms. Specifically, one-neighborhood strategy refers to moving only one vertex in the current partition set to another partition set to get a new set of partition results. However, the two-interchanges strategy refers to exchanging a vertex in the current partition set with a vertex in another partition set to obtain a new set of partition results.

On the one hand, we can use the two strategies of one-neighborhood and two-interchanges to adjust the infeasible vertices to feasible, and on the other hand, we also use them to improve the quality of the partition results. That is, for the $n_{P} \left(n_{P} \geq 2\right)$ feasible partitions $P_{1}, P_{2}, \ldots, P_{n_{P}} $ obtained by Algorithm 2, two partitions $P_{i}$ and $P_{j},$ where $i\in \left[n_{P}\right], j\in \left[n_{P}\right]\backslash i,$ are randomly selected from them, and the two strategies of one-neighborhood and two-interchanges are applied to adjust the vertices, and all feasible vertices are traversed to find the best partitioning result that can improve the objective function value after vertex tuning, provided that the vertex weight constraint is satisfied. Then, a heuristic algorithm to improve the quality of graph bipartition by using one-neighborhood and two-interchanges adjustment strategies is outlined in algorithm 3.

\vspace{0.6cm}
\begin{algorithm}
	\caption{Heuristic algorithm for improving partition quality}
\noindent
\textbf{Input:} Two randomly selected feasible partitions $P_{i}$ and $P_{j},$ $i\in \left[n_{P}\right], j\in \left[n_{P}\right]\backslash i;$\\
\textbf{Output:} Improved partitions $P^{*}_{i}$ and $P^{*}_{j},$ $i\in \left[n_{P}\right], j\in \left[n_{P}\right]\backslash i.$\\  \indent
\textbf{Step 0.} Calculate the objective function $\bar{f}\left({y}_{1}, {y}_{2}\right)$ corresponding to partitions $P_{i}$ and $P_{j}$ according to \\  \indent
\hspace*{1.2cm} model  \eqref{(e2.4)}, and set $\bar{P}_{i} := P_{i},$ $\bar{P}_{j} := P_{j}.$\\  \indent
\textbf{Step 1.} One-neighborhood adjustment to vertex $r:$\\  \indent
\hspace*{1.2cm} $P_{i} \leftarrow \bar{P}_{i} - \{r\};$ $P_{j} \leftarrow \bar{P}_{j} + \{r\},$  $i\in \left[n_{P}\right], j\in \left[n_{P}\right]\backslash i.$ \\  \indent
\textbf{Step 2.} Calculate the objective function ${f}\left({y}_{1}, {y}_{2}\right)$ corresponding to current  partitions $P_{i}$ and $P_{j}$ accord-\\  \indent
\hspace*{1.2cm} ing to model  \eqref{(e2.4)}.\\  \indent
\textbf{Step 3.} Judge whether the partition is feasible and whether the objective function decreases after \\  \indent
\hspace*{1.2cm} adjusting the vertices:\\  \indent
\hspace*{1.2cm} If ~~ ${y}_{1}^{T}{B} \leq {U}_1,$   ${y}_{2}^{T}{B} \leq {U}_2$ and ${f}\left({y}_{1}, {y}_{2}\right) < \bar{f}\left({y}_{1}, {y}_{2}\right)$ \\  \indent
\hspace*{1.8cm} Set $\bar{f}\left({y}_{1}, {y}_{2}\right) := {f}\left({y}_{1}, {y}_{2}\right),$ $\bar{P}_{i} := P_{i},$ $\bar{P}_{j} := P_{j};$ go to \textbf{Step 1.} \\  \indent
\hspace*{1.2cm} else \\  \indent
\hspace*{1.8cm} $\bar{P}_{i},$ $\bar{P}_{j}$ and $\bar{f}\left({y}_{1}, {y}_{2}\right)$ are given, go to \textbf{Step 4.} \\  \indent
\hspace*{1.2cm} end \\  \indent
\textbf{Step 4.} Two-interchanges adjustment to  vertices $s$ and $t:$ \\  \indent
\hspace*{1.2cm} $P_{i} \leftarrow \bar{P}_{i} - \{s\} + \{t\};$ $P_{j} \leftarrow \bar{P}_{j} + \{s\} - \{t\},$  $i\in \left[n_{P}\right], j\in \left[n_{P}\right]\backslash i.$ \\  \indent
\textbf{Step 5.} Calculate the objective function ${f}\left({y}_{1}, {y}_{2}\right)$ corresponding to current  partitions $P_{i}$ and $P_{j}$ accord-\\  \indent
\hspace*{1.2cm} ing to model  \eqref{(e2.4)}.\\  \indent
\textbf{Step 6.} Judge whether the partition is feasible and whether the objective function decreases after \\  \indent
\hspace*{1.2cm} adjusting the vertices:\\  \indent
\hspace*{1.2cm} If ~~ ${y}_{1}^{T}{B} \leq {U}_1,$   ${y}_{2}^{T}{B} \leq {U}_2$ and ${f}\left({y}_{1}, {y}_{2}\right) < \bar{f}\left({y}_{1}, {y}_{2}\right)$ \\  \indent
\hspace*{1.8cm} Set $\bar{f}\left({y}_{1}, {y}_{2}\right) := {f}\left({y}_{1}, {y}_{2}\right),$ $\bar{P}_{i} := P_{i},$ $\bar{P}_{j} := P_{j};$ go to \textbf{Step 4.} \\  \indent
\hspace*{1.2cm} else \\  \indent
\hspace*{1.8cm} Output $P^{*}_{i} := \bar{P}_{i}$ and $P^{*}_{j} := \bar{P}_{j},$ stop. \\  \indent
\hspace*{1.2cm} end \\  \indent
\end{algorithm}
\vspace{0.3cm}

\subsection{Recursive bipartition algorithm for multiple partitioning of graphs}
The multi-partitioning problem of a graph is to divide the vertices of the graph into multiple sets, and minimize the sum of the edge weights connecting the vertices of different sets under the condition that the sum of the $j$-th dimensional weights of the vertices of different sets does not exceed the resource limit.
In this paper, we use a recursive bipartition algorithm to solve the multi-partitioning problem of graphs. Each execution of the recursive bipartition algorithm means that the problem size of calling Algorithm 1 to solve the unconstrained model \eqref{(e2.11)} is much smaller than the problem solved in the previous, thus improving the efficiency of the graph partitioning problem. Algorithm 4 describes the flow of solving the graph partitioning problem with the recursive bipartition algorithm.

\vspace{0.6cm}
\begin{algorithm}
	\caption{Recursive bipartition algorithm for multiple partitioning of graphs (Rba)}
\noindent
\textbf{Step 0.} Call \textbf{Algorithm 1} to calculate the optimal solution of $\theta^{*}$ the problem \eqref{(e2.4)}. \\
\textbf{Step 1.} Call \textbf{Algorithm 2} to obtain the initial feasible partitions $P_{1}$ and $P_{2}$ of the graph bipartition \\
\hspace*{1.2cm} problem according to $\theta^{*}.$ \\
\textbf{Step 2.} Judge whether recursive bipartition is required:\\
\hspace*{1.2cm} If ~~ Resources are added in partition $P_{i} (i\in [2]),$ \\
\hspace*{1.8cm} \textbf{(1)} Call \textbf{Algorithm 2} to continue the bipartition for the vertices in partition $P_{i} (i\in [2]);$  \\
\hspace*{1.8cm} \textbf{(2)} Repeat the \textbf{Step 2.} until the $n_{P}$-partition result of the graph satisfying the condition \\
\hspace*{1.8cm} is obtained, i.e., $P_{1}, P_{2}, \ldots, P_{n_{P}} \left(n_{P}>2\right),$ go to \textbf{Step 3.}\\
\hspace*{1.2cm} else \\
\hspace*{1.8cm} Give the $n_{P}$-partition result of graph, i.e., $P_{1}, P_{2}, \ldots, P_{n_{P}} \left(n_{P}=2\right),$ go to \textbf{Step 3.}\\
\hspace*{1.2cm} end \\
\textbf{Step 3.} Call \textbf{Algorithm 3} to perform one-neighborhood adjustment and two-interchanges adjustment \\
\hspace*{1.2cm}for two randomly selected partitions $P_{i}$ and $P_{j},$ $i\in \left[n_{P}\right], j\in \left[n_{P}\right]\backslash i,$ and then get the improved\\
\hspace*{1.2cm} partitions $P^{*}_{i}$ and $P^{*}_{j}.$\\
\textbf{Step 4.} Output all improved partitions $P^{*}_{1}, P^{*}_{2}, \ldots, P^{*}_{n_{P}} .$
\end{algorithm}
\vspace{0.3cm}

\subsection{Time Complexity Analysis of the Algorithm}
This paper primarily utilizes the Subspace Minimization Conjugate Gradient (SMCG) method combined with the Recursive Bisection Algorithm to solve the graph partitioning problem. Therefore, the overall time complexity encompasses the time complexities of both algorithms.

The Recursive Bisection Algorithm divides the graph into $\log_2 n$ layers, with a recursion depth of $\log_2 n$, where $n$ represents the number of vertices in the graph. Thus, its time complexity is $O(\log n)$. At each recursion level, the Subspace Minimization Conjugate Gradient Method is invoked to solve a problem whose scale corresponds to the dimension of the current subgraph. The time complexity of the SMCG method is primarily influenced by the number of vertices in the current graph and the sparsity of the Laplacian matrix. Assuming the condition number $\kappa$ of the Laplacian matrix remains constant across all iterations, the time complexity of the SMCG method at each level is $O(m\sqrt{\kappa})$, where
$m$ represents the number of edges in the current subgraph. For sparse graphs, $m=O(n)$, while for dense graphs, $m=O(n^{2})$.  Consequently, the overall time complexity of the entire algorithm is $O(m\sqrt{\kappa}\log n )$.

\section{Numerical Experiments}\label{sec4}
This section is devoted to reporting some numerical results of the recursive bipartition Algorithm 4 (Rba), the Vc+2opt \cite{Wiegele2022} and  hypergraph partitioner \cite{Li2022}. To evaluate the performance of our algorithm, we conducted numerical experiments on two types of instances, including graph partition problem under knapsack constraints and hypergraph partitioning problem in very large scale integrated circuit design. All experiments are conducted on a desktop computer with an Intel(R) Core(TM) i5-6500 3.20GHz CPU and 16GB RAM.

\subsection{Graph Partition Problem Under Knapsack Constraints (GPKC)}
The first set of instances for GPKC problem are described in \cite{Nguyen2016} and the way we generate graph data is the same as that in reference \cite{Wiegele2022}, as follows:\\
1. Randomly select the edges of a complete graph with probabilities of 20\%, 50\% and 80\%, and note the corresponding instances as GPKCrand20, GPKCrand50 and GPKCrand80, respectively.\\
2. Non-zero edge weights are chosen as integers from the interval $(0,100].$\\
3. The vertex weights are one-dimensional and are chosen from integers in the interval $(0,1000].$\\
4. The number of vertices of the graph is 100, 200, 300, 400 and 500 respectively.

In the first set of numerical experiments, tables \ref{table1}, \ref{table2} and \ref{table3} give a detailed comparison of  edge cutting weight and CPU running time of different capacity bounds for the sparse graph  GPKCrand20, GPKCrand50 and GPKCrand80, respectively. The maximum number of iterations of ASMCG\_BB in algorithm Rba and the extended ADMM in algorithm Vc+2opt is set to 10000, the stopping tolerance is set to $10^{-5},$ and the equilibrium factor $\rho=5.$ In tables \ref{table1}, \ref{table2} and \ref{table3}, a $*$ indicates for the extended ADMM or ASMCG\_BB that the maximum number of iterations is reached. Cut-size denotes the total weight of the minimized cut edge obtained by the algorithm Rba  and the algorithm Vc+2opt, and the corresponding time is noted as CPU time . Gap represents the relative error of Cut-size solved by the two algorithms, i.e., Gap=(Rba-(Vc+2opt))/(Vc+2opt). In order to ensure the scientificity of the experimental results, this paper conducts 20 experiments on all the test cases, taking the minimum objective function value and the average CPU time.

As can be seen from tables \ref{table1}, \ref{table2} and \ref{table3}, compared with Vc+2opt, Rba only requires a very small amount of CPU time to obtain comparable partitioning results for any of GPKCrand20, GPKCrand50 and GPKCrand80, and the time required to partition the graph for Rba gradually increases as the capacity limit decreases. It is easy to see that if the gap is negative, it means that our method gets better partitioning results than Vc+2opt, otherwise, it means that our method is slightly inferior to Vc+2opt. At this time, the gap between the partitioning results obtained by the two algorithms can be measured by the size of Gaps. Table \ref{table1} shows that Gaps between Rba and Vc+2opt are less than 7.5\% for GPKCrand20, they are less than 3.5\% for GPKCrand50, see Table \ref{table2}, and for GPKCrand80, the gaps are less than 2\%, see Table \ref{table3}. As the results in these three tables show, Rba takes at most no more than 35 seconds and at least less than 1 second, while Vc+2opt takes at most more than an hour and at least 33 seconds. In addition, we note that it is more difficult for Vc+2opt to compute sparse instances because the maximum number of iterations is reached for GPKCrand20 much more often than for GPKCrand80 or GPKCrand50, but Rba performs better in this case, it means that ASMCG\_BB in Rba needs only a small number of iterations to obtain a better solution.

Figures 1-6 present the comparative line plots of Cut-size and CPU time between the two methods across different problem sizes and sparse graphs. Specifically, Figures 1, 3, and 5 illustrate the performance of the algorithm for Cut-size, while Figures 2, 4, and 6 depict the results for CPU time.
As shown in Figures 1, 3, and 5, the minimum cuts obtained by both Algorithm Rab and Algorithm VC+2opt increase with the problem size. However, Algorithm Rab consistently achieves smaller minimum cuts compared to Algorithm VC+2opt.
From Figures 2, 4, and 6, it is evident that for problems of the same scale, Algorithm Rab requires significantly less CPU time than Algorithm VC+2opt.
Overall, Figures 1-6 demonstrate that Algorithm Rab outperforms Algorithm VC+2opt in terms of both solution quality (Cut-size) and computational efficiency (CPU time).

\begin{table}[htp]
\caption{Numerical results for randomly generated graphs on GPKC problems GPKCrand20}\centering
\label{table1}
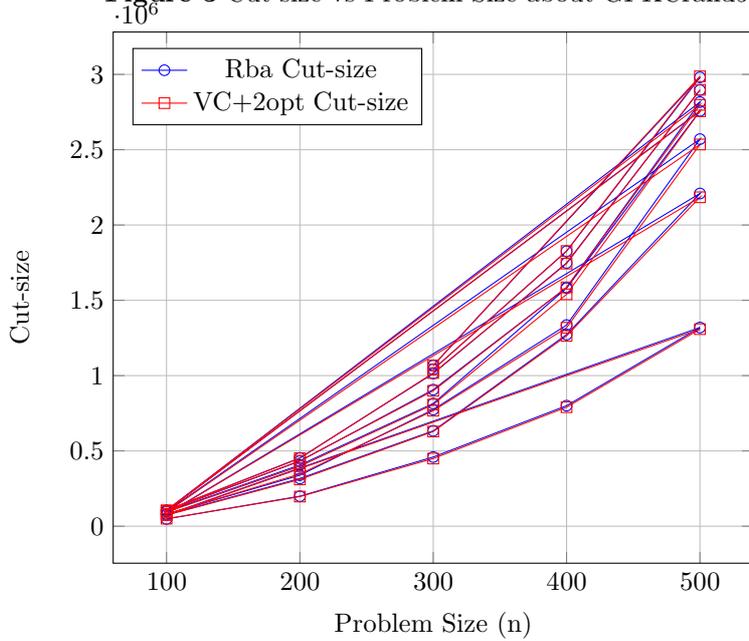
\begin{tabular}{@{}cccccccc@{}}
\toprule
\multirow{2}{*}{n} & \multirow{2}{*}{Capacity bound} & \multicolumn{2}{c}{Rba}        &  & \multicolumn{2}{c}{VC+2opt} & \multirow{2}{*}{Gap(\%)} \\ \cmidrule(lr){3-4} \cmidrule(lr){6-7}
                   &                                 & Cut-size          & CPU time(s) &  & Cut-size    & CPU time(s)    &                          \\ \midrule
100                & 26472                           & \textbf{16511}   & 0.026       &  & 16751*     & 169.979        & -1.4328                  \\
                   & 14310                           & \textbf{25680}   & 0.047       &  & 25764*     & 179.105        & -0.3260                  \\
                   & 11397                           & \textbf{28033}   & 0.092       &  & 28162*     & 180.041        & -0.4581                  \\
                   & 6349                            & \textbf{33407}   & 0.327       &  & 34228*     & 177.654        & -2.3986                  \\
                   & 3973                            & \textbf{36398}   & 0.458       &  & 37586*     & 177.295        & -3.1608                  \\
                   & 3082                            & \textbf{32007}   & 1.075       &  & 39187*     & 172.943        & -18.3224                 \\
200                & 49329                           & 78638            & 0.136       &  & 73315*     & 597.039        & 7.2605                   \\
                   & 26138                           & 115099           & 0.311       &  & 114481*    & 614.878        & 0.5398                   \\
                   & 21312                           & 126183           & 0.349       &  & 125430*    & 614.270        & 0.6003                   \\
                   & 11665                           & \textbf{146659}  & 0.774       &  & 148791*    & 609.839        & -1.4329                  \\
                   & 7007                            & \textbf{159034}  & 1.461       &  & 161382     & 150.229        & -1.4549                  \\
                   & 3747                            & \textbf{171603}  & 4.002       &  & 173093*    & 606.875        & -0.8608                  \\
300                & 73762                           & 183711           & 0.314       &  & 173857*    & 1312.709       & 5.6679                   \\
                   & 49938                           & 242663           & 0.423       &  & 237034*    & 1334.815       & 2.3748                   \\
                   & 31051                           & 299921           & 0.588       &  & 293800*    & 1359.697       & 2.0834                   \\
                   & 25979                           & 314105           & 0.706       &  & 310223*    & 1358.940       & 1.2514                   \\
                   & 16089                           & 349731           & 1.151       &  & 346213     & 1367.107       & 1.0161                   \\
                   & 6616                            & \textbf{387074}  & 4.380       &  & 389907     & 402.592        & -0.7266                  \\
                   & 4277                            & \textbf{401731}  & 7.582       &  & 402726     & 455.539        & -0.2471                  \\
                   & 2829                            & \textbf{410915}  & 19.638      &  & 411635     & 554.014        & -0.1749                  \\
400                & 96066                           & 322528           & 0.545       &  & 306186*    & 2431.096       & 5.3373                   \\
                   & 49622                           & 489681           & 0.721       &  & 484819*    & 2469.894       & 1.0028                   \\
                   & 40512                           & \textbf{528377}  & 0.952       &  & 530477*    & 2474.901       & -0.3959                  \\
                   & 26495                           & 588817           & 1.246       &  & 581763*    & 2467.382       & 1.2125                   \\
                   & 22338                           & 610502           & 1.564       &  & 608674*    & 2495.813       & 0.3003                   \\
                   & 11274                           & \textbf{671581}  & 2.639       &  & 674964*    & 2519.720       & -0.5012                  \\
                   & 7300                            & \textbf{700943}  & 6.913       &  & 705724     & 843.973        & -0.6775                  \\
500                & 133561                          & 495027           & 0.886       &  & 490984*    & 4040.566       & 0.8234                   \\
                   & 54084                           & 862724           & 1.247       &  & 848777     & 1784.043       & 1.6432                   \\
                   & 28866                           & 1014656          & 1.978       &  & 982437     & 1646.458       & 3.2795                   \\
                   & 14840                           & \textbf{1072500} & 4.674       &  & 1082707    & 1463.913       & -0.9427                  \\
                   & 12923                           & \textbf{1087406} & 5.319       &  & 1092863    & 1529.387       & -0.4993                  \\
                   & 7468                            & \textbf{1129457} & 12.757      &  & 1132038    & 1611.319       & -0.2280                  \\
                   & 3975                            & \textbf{1166265} & 33.131      &  & 1170041    & 2006.805       & -0.3227                  \\ \bottomrule
\end{tabular}
\end{table}

\begin{tikzpicture}
\begin{axis}[
    xlabel={Problem Size (n)},
    ylabel={Cut-size},
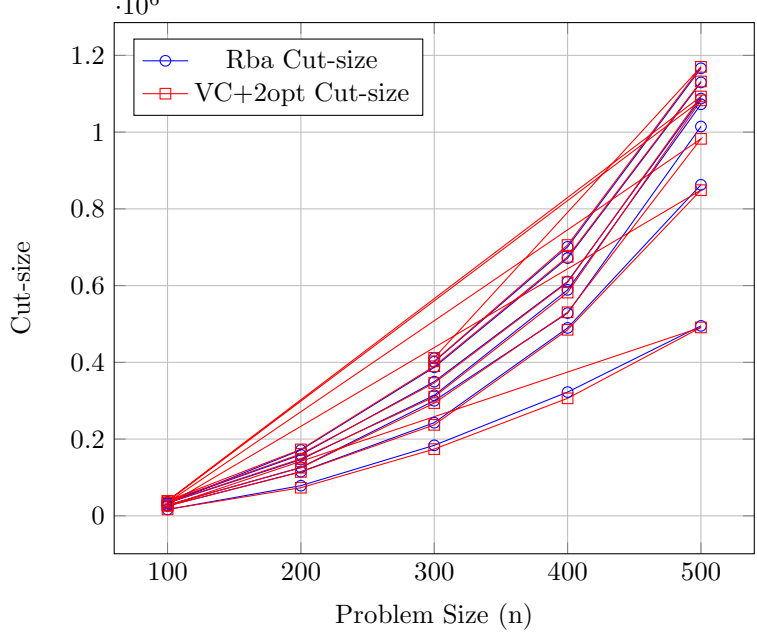
    title={\textbf{Figure 1} Cut-size vs Problem Size about GPKCrand20},
    legend pos=north west,
    grid=major,
    xtick={100,200,300,400,500},
    ymajorgrids=true,
    xmajorgrids=true
]
\addplot[blue, mark=o] coordinates {

(100, 16511) (200, 78638)(300, 183711)(400, 322528)(500, 495027)

(100, 25680)  (200, 115099) (300, 242663)(400, 489681)(500, 862724)

(100, 28033) (200, 126183)(300, 299921)(400, 528377)(500, 1014656)

(100, 33407) (200, 146659)(300, 314105)(400, 588817)(500, 1072500)

(100, 36398) (200, 159034)(300, 349731)(400, 610502)(500, 1087406)

(100, 32007)(200, 171603)(300, 387074)(400, 671581)(500, 1129457)

(300, 401731) (400, 700943)(500, 1166265)

(300, 410915)

};
\addplot[red, mark=square] coordinates {
(100, 16751) (200, 73315)(300, 173857)(400, 306186)(500, 490984)
(100, 25764) (200, 114481)(300, 237034) (400, 484819) (500, 848777)
(100, 28162) (200, 125430)(300, 293800)(400, 530477) (500, 982437)
(100, 34228) (200, 148791)(300, 310223)(400, 581763)(500, 1082707)
(100, 37586) (200, 161382)(300, 346213)(400, 608674)(500, 1092863)
(100, 39187) (200, 173093)(300, 389907)(400, 674964)  (500, 1132038)

(300, 402726) (400, 705724) (500, 1170041)
(300, 411635)

};
\legend{Rba Cut-size, VC+2opt Cut-size}
\end{axis}
\end{tikzpicture}

\begin{tikzpicture}
\begin{axis}[
    xlabel={Problem Size (n)},
    ylabel={CPU Time (s)},
    title={\textbf{Figure 2 }CPU Time vs Problem Size about GPKCrand20(Log Scale)},
    legend pos=north west,
    grid=major,
    xtick={100,200,300,400,500},
    ymajorgrids=true,
    xmajorgrids=true,
    ymode=log, 
    log ticks with fixed point
]
\addplot[blue, mark=o] coordinates {

(100, 0.026) (200, 0.136)(300, 0.314)(400, 0.545)(500, 0.886)
(100, 0.047) (200, 0.311)(300, 0.423)(400, 0.721) (500, 1.247)
(100, 0.092) (200, 0.349)(300, 0.588)(400, 0.952)(500, 1.978)
(100, 0.327) (200, 0.774) (300, 0.706)(400, 1.246) (500, 4.674)
(100, 0.458) (200, 1.461)(300, 1.151)(400, 1.564) (500, 5.319)
(100, 1.075)(200, 4.002)(300, 4.380)(400, 2.639)  (500, 12.757)

(300, 7.582)(400, 6.913)(500, 33.131)

(300, 19.638)

};
\addplot[red, mark=square] coordinates {

(100, 169.979) (200, 597.039)(300, 1312.709)(400, 2431.096)(500, 4040.566)
(100, 179.105) (200, 614.878)(300, 1334.815) (400, 2469.894) (500, 1784.043)
(100, 180.041) (200, 614.270)(300, 1359.697) (400, 2474.901)(500, 1646.458)
(100, 177.654) (200, 609.839)(300, 1358.940)(400, 2467.382) (500, 1463.913)
(100, 177.295) (200, 150.229)(300, 1367.107) (400, 2495.813) (500, 1529.387)
(100, 172.943)(200, 606.875)(300, 402.592)(400, 2519.720)(500, 1611.319)

(300, 455.539)(400, 843.973)(500, 2006.805)
(300, 554.014)

};
\legend{Rba CPU Time, VC+2opt CPU Time}
\end{axis}
\end{tikzpicture}

\begin{table}[htp]
\caption{Numerical results for randomly generated graphs on GPKC problems GPKCrand50}\centering
\label{table2}
\begin{tabular}{@{}cccccccc@{}}
\toprule
\multirow{2}{*}{n} & \multirow{2}{*}{Capacity bound} & \multicolumn{2}{c}{Rba}        &  & \multicolumn{2}{c}{VC+2opt} & \multirow{2}{*}{Gap(\%)} \\ \cmidrule(lr){3-4} \cmidrule(lr){6-7}
                   &                                 & Cut-size          & CPU time(s) &  & Cut-size     & CPU time(s)   &                          \\ \midrule
100                & 27572                           & \textbf{48379}   & 0.032       &  & 49743*      & 169.292       & -2.7421                  \\
                   & 14639                           & 75405            & 0.090       &  & 73521*      & 170.525       & 2.5625                   \\
                   & 11686                           & 80950            & 0.118       &  & 80192       & 119.193       & 0.9452                   \\
                   & 6526                            & 97039            & 0.219       &  & 93793       & 36.606        & 3.4608                   \\
                   & 4040                            & 102615           & 0.461       &  & 102034      & 38.212        & 0.5694                   \\
                   & 3179                            & \textbf{105573}  & 0.593       &  & 105737      & 42.276        & -0.1551                  \\
200                & 48935                           & 199268           & 0.164       &  & 196678*     & 602.008       & 1.3169                   \\
                   & 25685                           & 316844           & 0.253       &  & 310687*     & 618.162       & 1.9817                   \\
                   & 21504                           & \textbf{338620}           & 0.338       &  & 346751*     & 613.326       & -2.3449                  \\
                   & 10742                           & 406369           & 0.819       &  & 398454      & 143.755       & 1.9864                   \\
                   & 6034                            & \textbf{432589}  & 4.538       &  & 434442*     & 609.051       & -0.4265                  \\
                   & 3940                            & 453419           & 3.571       &  & 451562*     & 602.538       & 0.4112                   \\
300                & 76829                           & 459421           & 0.246       &  & 449299*     & 1305.923      & 2.2528                   \\
                   & 49664                           & 633274           & 0.408       &  & 628741      & 432.045       & 0.7210                   \\
                   & 32370                           & 774369           & 0.580       &  & 767085      & 421.686       & 0.9496                   \\
                   & 27123                           & 814752           & 0.785       &  & 807901      & 412.435       & 0.8480                   \\
                   & 16670                           & 906547           & 1.217       &  & 899218      & 388.895       & 0.8150                   \\
                   & 6371                            & \textbf{1013686} & 5.384       &  & 1017229     & 421.242       & -0.3483                  \\
                   & 4428                            & 1041398          & 7.154       &  & 1041040     & 450.359       & 0.0344                   \\
                   & 2649                            & \textbf{1067053} & 17.813      &  & 1067431     & 550.762       & -0.0354                  \\
400                & 97431                           & 800315           & 0.494       &  & 790271*     & 2415.262      & 1.2710                   \\
                   & 48514                           & 1272445          & 0.844       &  & 1263905     & 966.465       & 0.6757                   \\
                   & 44456                           & 1335383          & 0.830       &  & 1317855     & 954.913       & 1.3300                   \\
                   & 25558                           & 1580893          & 1.200       &  & 1539624     & 910.983       & 2.6805                   \\
                   & 22264                           & \textbf{1585208} & 1.626       &  & 1585733     & 920.929       & -0.0331                  \\
                   & 11477                           & 1748040          & 3.850       &  & 1743822     & 934.501       & 0.2419                   \\
                   & 6800                            & \textbf{1825977} & 7.551       &  & 1827386     & 988.103       & -0.0771                  \\
500                & 131286                          & 1319301          & 0.796       &  & 1308587*    & 3943.253      & 0.8187                   \\
                   & 56101                           & 2209206          & 1.190       &  & 2184011     & 1662.666      & 1.1536                   \\
                   & 29244                           & 2570289          & 2.288       &  & 2534630     & 1555.525      & 1.4069                   \\
                   & 15747                           & 2757429          & 5.142       &  & 2755417     & 1592.084      & 0.0730                   \\
                   & 13176                           & 2816204          & 5.485       &  & 2797682     & 1587.139      & 0.6620                   \\
                   & 7920                            & \textbf{2895821} & 12.496      &  & 2897871     & 1706.526      & -0.0707                  \\
                   & 4249                            & \textbf{2980653} & 29.527      &  & 2987303     & 1947.095      & -0.2226                  \\ \bottomrule
\end{tabular}
\end{table}

\begin{tikzpicture}
\begin{axis}[
    xlabel={Problem Size (n)},
    ylabel={Cut-size},
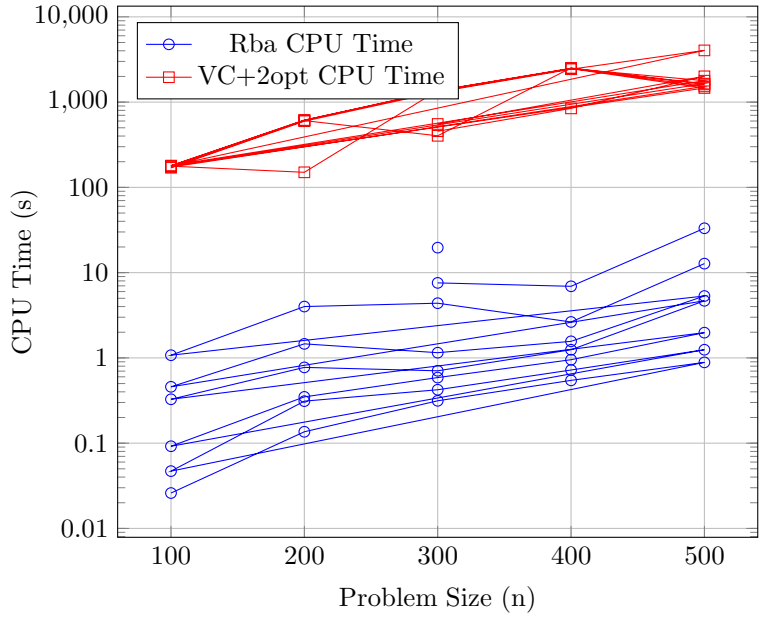
    title={\textbf{Figure 3} Cut-size vs Problem Size about GPKCrand50},
    legend pos=north west,
    grid=major,
    xtick={100,200,300,400,500},
    ymajorgrids=true,
    xmajorgrids=true
]
\addplot[blue, mark=o] coordinates {

(100, 48379) (200, 199268)(300, 459421)(400, 800315) (500, 1319301)
(100, 75405) (200, 316844)(300, 633274)(400, 1272445)(500, 2209206)
(100, 80950) (200, 338620)(300, 774369) (400, 1335383)(500, 2570289)
(100, 97039) (200, 406369)(300, 814752)(400, 1580893)(500, 2757429)
(100, 102615) (200, 432589)(300, 906547)(400, 1585208)(500, 2816204)
(100, 105573) (200, 453419)(300, 1013686)(400, 1748040)(500, 2895821)

(300, 1041398) (400, 1825977) (500, 2980653)
(300, 1067053)

};
\addplot[red, mark=square] coordinates {
(100, 49743) (200, 196678)(300, 449299)(400, 790271)(500, 1308587)
(100, 73521) (200, 310687)(300, 628741)(400, 1263905)(500, 2184011)
(100, 80192) (200, 346751)(300, 767085)(400, 1317855)(500, 2534630)
(100, 93793) (200, 398454)(300, 807901)(400, 1539624)(500, 2755417)
(100, 102034) (200, 434442)(300, 899218) (400, 1585733) (500, 2797682)
(100, 105737)(200, 451562)(300, 1017229)(400, 1743822) (500, 2897871)

(300, 1041040)(400, 1827386)(500, 2987303)
(300, 1067431)

};
\legend{Rba Cut-size, VC+2opt Cut-size}
\end{axis}
\end{tikzpicture}

\begin{tikzpicture}
\begin{axis}[
    xlabel={Problem Size (n)},
    ylabel={CPU Time (s)},
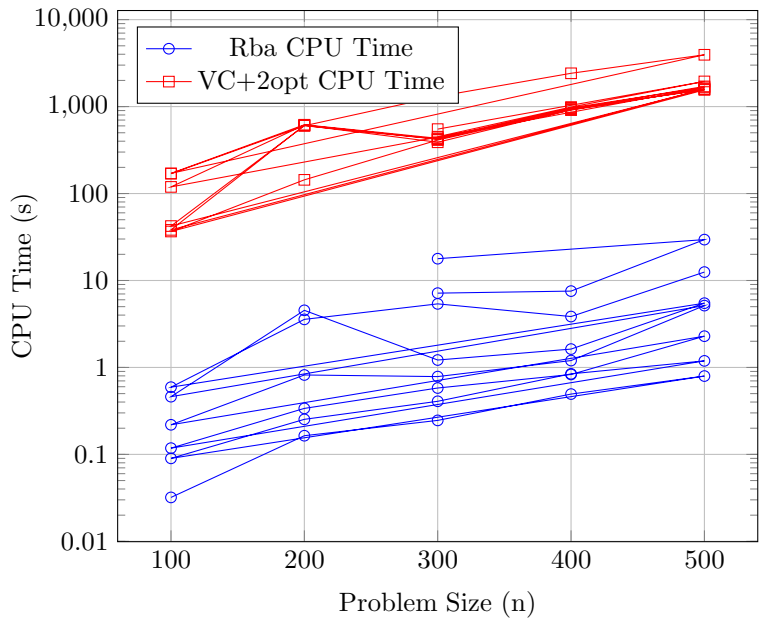
    title={\textbf{Figure 4 }CPU Time vs Problem Size about GPKCrand50(Log Scale)},
    legend pos=north west,
    grid=major,
    xtick={100,200,300,400,500},
    ymajorgrids=true,
    xmajorgrids=true,
    ymode=log, 
    log ticks with fixed point
]
\addplot[blue, mark=o] coordinates {
(100, 0.032) (200, 0.164)(300, 0.246)(400, 0.494)(500, 0.796)
(100, 0.090) (200, 0.253)(300, 0.408)(400, 0.844)(500, 1.190)
(100, 0.118) (200, 0.338)(300, 0.580)(400, 0.830)(500, 2.288)
(100, 0.219) (200, 0.819)(300, 0.785)(400, 1.200)(500, 5.142)
(100, 0.461) (200, 4.538)(300, 1.217)(400, 1.626)(500, 5.485)
(100, 0.593) (200, 3.571)(300, 5.384)(400, 3.850)(500, 12.496)

(300, 7.154) (400, 7.551)(500, 29.527)
(300, 17.813)

};
\addplot[red, mark=square] coordinates {

(100, 169.292) (200, 602.008)(300, 1305.923)(400, 2415.262)(500, 3943.253)
(100, 170.525) (200, 618.162)(300, 432.045)(400, 966.465)(500, 1662.666)
(100, 119.193) (200, 613.326) (300, 421.686)(400, 954.913)(500, 1555.525)
(100, 36.606) (200, 143.755)(300, 412.435)(400, 910.983)(500, 1592.084)
(100, 38.212) (200, 609.051)(300, 388.895)(400, 920.929)(500, 1587.139)
(100, 42.276)(200, 602.538)(300, 421.242)(400, 934.501)(500, 1706.526)

     (300, 450.359)(400, 988.103) (500, 1947.095)
      (300, 550.762)

};
\legend{Rba CPU Time, VC+2opt CPU Time}
\end{axis}
\end{tikzpicture}

\begin{table}[htp]
\caption{Numerical results for randomly generated graphs on GPKC problems GPKCrand80}\centering
\label{table3}
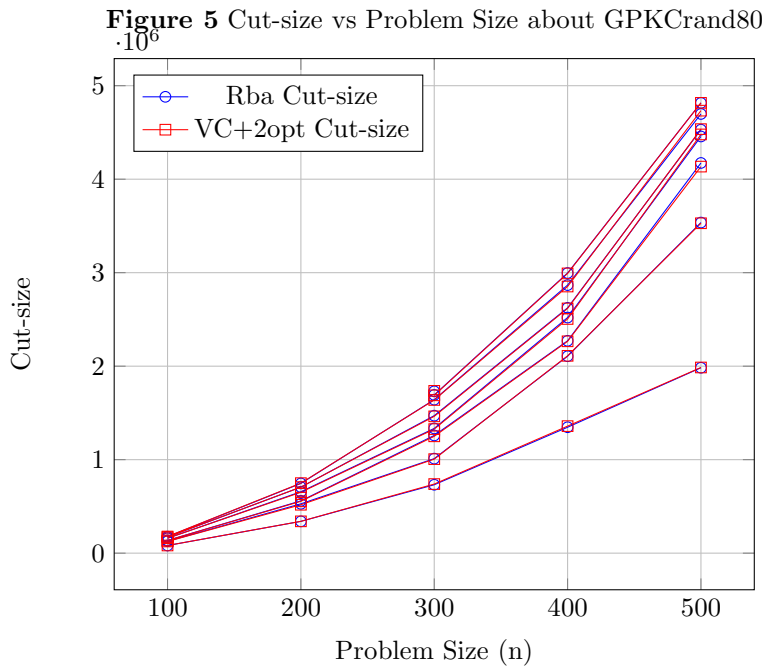
\begin{tabular}{@{}cccccccc@{}}
\toprule
\multirow{2}{*}{n} & \multirow{2}{*}{Capacity bound} & \multicolumn{2}{c}{Rba}        &  & \multicolumn{2}{c}{VC+2opt} & \multirow{2}{*}{Gap(\%)} \\ \cmidrule(lr){3-4} \cmidrule(lr){6-7}
                   &                                 & Cut-size          & CPU time(s) &  & Cut-size     & CPU time(s)   &                          \\ \midrule
100                & 26536                           & \textbf{81897}   & 0.033       &  & 81979*      & 164.792       & -0.1000                  \\
                   & 14023                           & 127081           & 0.079       &  & 124979      & 33.837        & 1.6819                   \\
                   & 11434                           & 136103           & 0.118       &  & 135996      & 35.100        & 0.0787                   \\
                   & 6321                            & 159436           & 0.229       &  & 159434      & 37.565        & 0.0013                   \\
                   & 3668                            & \textbf{172700}  & 0.526       &  & 172867      & 43.821        & -0.0969                  \\
                   & 3043                            & \textbf{176601}  & 0.629       &  & 178090      & 52.689        & -0.8364                  \\
200                & 50084                           & \textbf{340066}  & 0.182       &  & 340370      & 265.389       & -0.0893                  \\
                   & 25787                           & 527658           & 0.252       &  & 517839      & 185.867       & 1.8961                   \\
                   & 21741                           & 559424           & 0.336       &  & 554627      & 193.917       & 0.8649                   \\
                   & 11404                           & 656842           & 0.687       &  & 651318      & 212.021       & 0.8481                   \\
                   & 6686                            & \textbf{706304}  & 1.867       &  & 706648      & 203.074       & -0.0487                  \\
                   & 3419                            & \textbf{745638}  & 6.814       &  & 748894      & 201.155       & -0.4348                  \\
300                & 73485                           & \textbf{731591}  & 0.256       &  & 739066*     & 1318.380      & -1.0114                  \\
                   & 49920                           & 1011057          & 0.342       &  & 1005754     & 635.548       & 0.5273                   \\
                   & 31134                           & 1258954          & 0.546       &  & 1249668     & 766.800       & 0.7431                   \\
                   & 25752                           & 1334236          & 0.657       &  & 1326910     & 766.566       & 0.5521                   \\
                   & 16767                           & 1470714          & 1.051       &  & 1463377     & 684.852       & 0.5014                   \\
                   & 6647                            & 1641797          & 3.810       &  & 1638592     & 513.306       & 0.1956                   \\
                   & 4210                            & 1692267          & 7.977       &  & 1688172     & 618.567       & 0.2426                   \\
                   & 2454                            & \textbf{1730890} & 21.228      &  & 1736053     & 692.605       & -0.2974                  \\
400                & 99348                           & \textbf{1346717} & 0.570       &  & 1358776*    & 2464.254      & -0.8875                  \\
                   & 51033                           & \textbf{2106341} & 0.729       &  & 2110561     & 1469.015      & -0.1999                  \\
                   & 41130                           & 2271413          & 0.881       &  & 2266076     & 1530.471      & 0.2355                   \\
                   & 28718                           & 2520045          & 1.314       &  & 2504817     & 1935.380      & 0.6079                   \\
                   & 22740                           & 2623899          & 1.728       &  & 2617677     & 1729.365      & 0.2377                   \\
                   & 11772                           & 2867205          & 3.776       &  & 2852005     & 1453.797      & 0.5330                   \\
                   & 6132                            & 2996214          & 10.323      &  & 2992314     & 1681.777      & 0.1303                   \\
500                & 132135                          & \textbf{1983188} & 0.694       &  & 1986063*    & 3933.542      & -0.1448                  \\
                   & 53186                           & 3536128          & 1.216       &  & 3528720*    & 4000.427      & 0.2099                   \\
                   & 26965                           & 4173638          & 2.081       &  & 4134651     & 3453.792      & 0.9429                   \\
                   & 15049                           & \textbf{4457006} & 5.357       &  & 4478716     & 3603.654      & -0.4847                  \\
                   & 12806                           & \textbf{4530025} & 6.377       &  & 4535627     & 3929.182      & -0.1235                  \\
                   & 7312                            & \textbf{4697683} & 14.474      &  & 4733054     & 2747.474      & -0.7473                  \\
                   & 4071                            & \textbf{4815995} & 34.057      &  & 4816388     & 2667.188      & -0.0082                  \\ \bottomrule
\end{tabular}
\end{table}

\begin{tikzpicture}
\begin{axis}[
    xlabel={Problem Size (n)},
    ylabel={Cut-size},
    title={\textbf{Figure 5} Cut-size vs Problem Size about GPKCrand80},
    legend pos=north west,
    grid=major,
    xtick={100,200,300,400,500},
    ymajorgrids=true,
    xmajorgrids=true
]
\addplot[blue, mark=o] coordinates {
    (100, 81897)
    (200, 340066)
    (300, 731591)
    (400, 1346717)
    (500, 1983188)

    (100, 127081)
    (200, 527658)
    (300, 1011057)
    (400, 2106341)
    (500, 3536128)

    (100, 136103)
    (200, 559424)
    (300, 1258954)
    (400, 2271413)
    (500, 4173638)

    (100, 159436)
    (200, 656842)
    (300, 1334236)
    (400, 2520045)
    (500, 4457006)

    (100, 172700)
    (200, 706304)
    (300, 1470714)
    (400, 2623899)
    (500, 4530025)

    (100, 176601)
    (200, 745638)
    (300, 1641797)
    (400, 2867205)
    (500, 4697683)

    (300, 1692267)
    (400, 2996214)
    (500, 4815995)

    (300, 1730890)

};
\addplot[red, mark=square] coordinates {
    (100, 81979)
    (200, 340370)
    (300, 739066)
    (400, 1358776)
    (500, 1986063)

    (100, 124979)
    (200, 517839)
    (300, 1005754)
    (400, 2110561)
    (500, 3528720)

    (100, 135996)
    (200, 554627)
    (300, 1249668)
    (400, 2266076)
    (500, 4134651)

    (100, 159434)
    (200, 651318)
    (300, 1326910)
    (400, 2504817)
    (500, 4478716)

    (100, 172867)
    (200, 706648)
    (300, 1463377)
    (400, 2617677)
    (500, 4535627)

    (100, 178090)
    (200, 748894)
    (300, 1638592)
    (400, 2852005)
    (500, 4733054)

    (300, 1688172)
    (400, 2992314)
    (500, 4816388)

    (300, 1736053)

};
\legend{Rba Cut-size, VC+2opt Cut-size}
\end{axis}
\end{tikzpicture}

\begin{tikzpicture}
\begin{axis}[
    xlabel={Problem Size (n)},
    ylabel={CPU Time (s)},
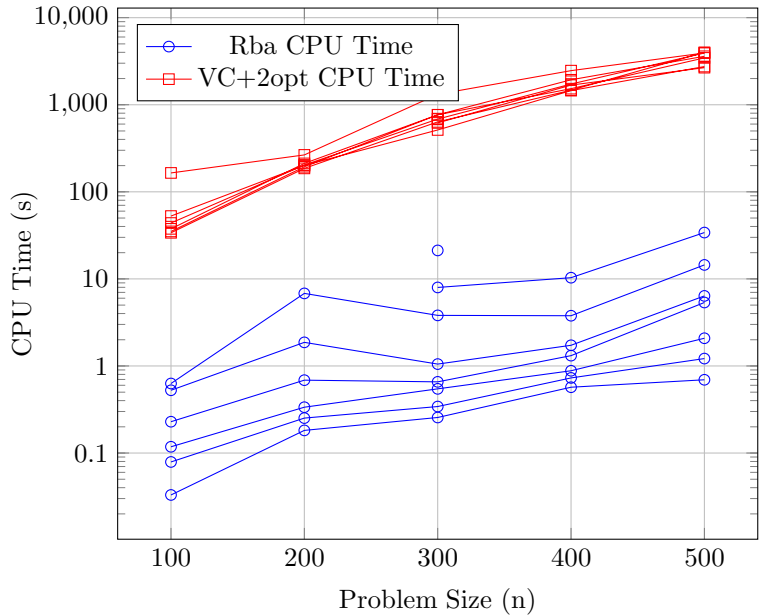
    title={\textbf{Figure 6 }CPU Time vs Problem Size about GPKCrand80(Log Scale)},
    legend pos=north west,
    grid=major,
    xtick={100,200,300,400,500},
    ymajorgrids=true,
    xmajorgrids=true,
    ymode=log, 
    log ticks with fixed point
]
\addplot[blue, mark=o] coordinates {
    (100, 0.033)
    (200, 0.182)
    (300, 0.256)
    (400, 0.570)
    (500, 0.694)

    (100, 0.079)
    (200, 0.252)
    (300, 0.342)
    (400, 0.729)
    (500, 1.216)

    (100, 0.118)
    (200, 0.336)
    (300, 0.546)
    (400, 0.881)
    (500, 2.081)

    (100, 0.229)
    (200, 0.687)
    (300, 0.657)
    (400, 1.314)
    (500, 5.357)

    (100, 0.526)
    (200, 1.867)
    (300, 1.051)
    (400, 1.728)
    (500, 6.377)

    (100, 0.629)
    (200, 6.814)
    (300, 3.810)
    (400, 3.776)
    (500, 14.474)

    (300, 7.977)
    (400, 10.323)
    (500, 34.057)

    (300, 21.228)
};
\addplot[red, mark=square] coordinates {
    (100, 164.792)
    (200, 265.389)
    (300, 1318.380)
    (400, 2464.254)
    (500, 3933.542)

    (100, 33.837)
    (200, 185.867)
    (300, 635.548)
    (400, 1469.015)
    (500, 4000.427)

    (100, 35.100)
    (200, 193.917)
    (300, 766.800)
    (400, 1530.471)
    (500, 3453.792)

    (100, 37.565)
    (200, 212.021)
    (300, 766.566)
    (400, 1935.380)
    (500, 3603.654)

    (100, 43.821)
    (200, 203.074)
    (300, 684.852)
    (400, 1729.365)
    (500, 3929.182)

    (100, 52.689)
    (200, 201.155)
    (300, 513.306)
    (400, 1453.797)
    (500, 2747.474)

    (300, 618.567)
    (400, 1681.777)
    (500, 2667.188)

    (300, 692.605)
};
\legend{Rba CPU Time, VC+2opt CPU Time}
\end{axis}
\end{tikzpicture}

\subsection{Hypergraph Partitioning in Very Large Scale Integrated Circuit Design (VLSI)}
The second set of instances are the circuit partitioning problem in the design of very large scale integrated circuits, derived from actual industrial design data. The circuit partitioning problem in VLSI is a typical hypergraph partitioning problem, where the nodes in the hypergraph represent the modules in the circuit and the hyperedges represent the signal networks in the circuit. Set $\bar{E}$ is the hypergraph edge, $\left|\bar{E}\right|$ denotes the number of vertices connected by the hypergraph edge, and $W_{\bar{E}}$ denotes the weight value of the hypergraph edge $\bar{E}$. If we decompose a signal network into $\left|\bar{E}\right| \left(\left|\bar{E}\right|-1\right)/2$ edges, and the weight value of each edge is assigned according to formula  $ W_{\bar{E}} / \left(\left|\bar{E}\right|-1\right)$, then the hypergraph problem can be transformed into a graph problem. The vertex weights of such problems are often multidimensional and often come with fixed vertex constraints.

This set of experiments is mainly to further compare the performance of algorithm Rba and hypergraph Partitioner proposed in reference \cite{Li2022} by testing the graph partitioning problem obtained from the transformation of hypergraph partitioning. Table \ref{table4} shows the numerical results of two methods for solving actual industrial design problems. Cases industry 1 - industry 6 are problems with fixed vertex constraints, where the fixed vertices are distributed in two different partitioned sets of $P_{1}$ and $P_{2}$, cases industry 7 and industry 8 have no fixed vertex constraints.  Cut-size denotes the total weight of the minimized cut edge obtained by the algorithm Rba and Partitioner \cite{Li2022}, and Gap represents the relative error of Cut-size obtained by the algorithm Rba and Partitioner\cite{Li2022}, i.e., Gap=(Rba-Partitioner)/Partitioner, and a `` - " indicates no fixed vertices or no relative error Gap is obtained.
As can be seen from Table \ref{table4}, Cut-size obtained by our method Rba is significantly better than that of Partitioner \cite{Li2022}, which indicates that our method is able to solve the hypergraph problem in circuit partitioning.

\begin{table}[htp]
\caption{Numerical results for Hypergraph partitioning problem}\centering
\label{table4}
\begin{tabular}{@{}ccccccccc@{}}
\toprule
\multicolumn{1}{c}{\multirow{2}{*}{Case}} & \multicolumn{1}{c}{\multirow{2}{*}{Nodes}} & \multicolumn{1}{c}{\multirow{2}{*}{Nets}} & \multicolumn{2}{c}{Fixed nodes}                  & \multirow{2}{*}{} & \multicolumn{2}{c}{Cut-size}     & \multirow{2}{*}{Gap(\%)}  \\ \cmidrule(lr){4-5} \cmidrule(lr){7-8}
\multicolumn{1}{c}{}                      & \multicolumn{1}{c}{}                       & \multicolumn{1}{c}{}                      & \multicolumn{1}{c}{$P_{1}$} & \multicolumn{1}{c}{$P_{2}$} &                   & \multicolumn{1}{c}{Rba}           & \multicolumn{1}{c}{Partitioner}      &          \\ \midrule
industry 1                                 & 53                                         & 40                                        & g4                     & g5                     &                   & 2             & 2                & 0.0000   \\
industry 2                                 & 274                                        & 9320                                      & g190                   & g93                    &                   & 1363          & 1363             & 0.0000   \\
industry 3                                 & 995                                        & 982                                       & g37                    & g51                    &                   & 30            & 30               & 0.0000   \\
industry 4                                 & 1173                                       & 281                                       & g30,g124               & g19                    &                   & \textbf{715}  & 748              & -4.4118  \\
industry 5                                 & 4773                                       & 22589                                     & g257                   & g4568                  &                   & 35            & 35               & 0.0000   \\
industry 6                                 & 11531                                      & 32665                                     & g3,g25,g227,g228       & g4                     &                   & \textbf{3843} & 3885             & -1.0811  \\
industry 7                                 & 19601                                      & 19583                                     & -                      & -                      &                   & 2067          & 1957(infeasible) & -\\
industry 8                                 & 69429                                      & 75195                                     & -                      & -                      &                   & \textbf{1920} & 11983            & -83.9773 \\ \bottomrule
\end{tabular}
\end{table}

In order to demonstrate that the proposed heuristic algorithm can indeed improve partition results, we test the algorithm with the data set of ISPD98 benchmark.
ISPD98 is currently used by researchers in academia and industry as a benchmark tool for circuit partitioning evaluation. The solution of circuit partition problem in ISPD9 benchmark can be converted to hypergraph partition problem. KaHyPar\cite{Akhremtsev17}  is one of the most efficient methods for solving hypergraph partitioning problems. Here, the result of KaHyPar solving ISPD98 benchmark is taken as the initial partition result, and then heuristic algorithm 3 is further used to improve it, so as to verify the effectiveness of the proposed heuristic algorithm.
Table \ref{table5} list the characteristic parameters after converting ISPD98 into hypergraph and result comparison data of the test data set.
In Table \ref{table5}, ``KaHyPar Cut-size"  represents the total weight of the minimized cut edge obtained by KaHyPar solving 18 problems in ISPD98  benchmark respectively. ``Improved  Cut-size" represents the result improved by algorithm 3 on the Cut-size in KaHyPar. It can be seen from Table \ref{table5} that the heuristic vertex adjustment strategy designed in this paper can improve most of the partition results of KaHyPar, which shows that the vertex adjustment strategy in the proposed algorithm 3 is effective.

\begin{table}[htp]
\caption{Improved comparison of results for KaHyPar partition}\centering
\label{table5}
\begin{tabular}{@{}cccccccc@{}}
\toprule
Circuit & Vertices & Hyperedges & Pins   & \begin{tabular}[c]{@{}c@{}}Vertex \\ Weights\end{tabular} & \begin{tabular}[c]{@{}c@{}}Hyperedge\\ Weights\end{tabular} & \begin{tabular}[c]{@{}c@{}}KaHyPar\\ Cut-size\end{tabular} & \begin{tabular}[c]{@{}c@{}}Improved\\ Cut-size\end{tabular} \\ \midrule
ibm01   & 12752    & 14111      & 50566  & 1                                                         & 1                                                           & 180                                                        & 180                                                         \\
ibm02   & 19601    & 19584      & 81199  & 1                                                         & 1                                                           & 262                                                        & 262                                                         \\
ibm03   & 23136    & 27401      & 93573  & 1                                                         & 1                                                           & 2689                                                       & \textbf{2587}                                               \\
ibm04   & 27507    & 31970      & 105859 & 1                                                         & 1                                                           & 2149                                                       & \textbf{2138}                                               \\
ibm05   & 29347    & 28446      & 126308 & 1                                                         & 1                                                           & 4474                                                       & \textbf{4422}                                               \\
ibm06   & 32498    & 34826      & 128182 & 1                                                         & 1                                                           & 2358                                                       & \textbf{2352}                                               \\
ibm07   & 45926    & 48117      & 175639 & 1                                                         & 1                                                           & 4479                                                       & \textbf{4268}                                               \\
ibm08   & 51309    & 50513      & 204890 & 1                                                         & 1                                                           & 5099                                                       & 5099                                                        \\
ibm09   & 53395    & 60902      & 222088 & 1                                                         & 1                                                           & 3179                                                       & \textbf{3055}                                               \\
ibm10   & 69429    & 75196      & 297567 & 1                                                         & 1                                                           & 4706                                                       & \textbf{4450}                                               \\
ibm11   & 70558    & 81454      & 280786 & 1                                                         & 1                                                           & 6161                                                       & \textbf{5832}                                               \\
ibm12   & 71076    & 77240      & 317760 & 1                                                         & 1                                                           & 8701                                                       & \textbf{8240}                                               \\
ibm13   & 84199    & 99666      & 357075 & 1                                                         & 1                                                           & 5283                                                       & \textbf{5036}                                               \\
ibm14   & 147605   & 152772     & 546816 & 1                                                         & 1                                                           & 9459                                                       & \textbf{8833}                                               \\
ibm15   & 161570   & 186608     & 715823 & 1                                                         & 1                                                           & 27821                                                      & \textbf{25497}                                              \\
ibm16   & 183484   & 190048     & 778823 & 1                                                         & 1                                                           & 31884                                                      & \textbf{28786}                                              \\
ibm17   & 185495   & 189581     & 860036 & 1                                                         & 1                                                           & 59645                                                      & \textbf{53996}                                              \\
ibm18   & 210613   & 201920     & 819697 & 1                                                         & 1                                                           & 44544                                                      & \textbf{38543}                                              \\ \bottomrule
\end{tabular}
\end{table}

From the results of the above  numerical experiments, apparently the proposed method is still quite effective for the graph partitioning problem. In addition, the method proposed in this paper is also very effective for large-scale graph partitioning with more than 50,000 vertices.

\section{Conclusions and Future Work}\label{sec5}
In this paper, a recursive bipartition  algorithm for solving graph partitioning problems with multidimensional vertex weight constraints and fixed vertex constraints is proposed based on the subspace minimization conjugate gradient algorithm. Firstly, for the 0-1 integer programming model of the graph partitioning problem, the complex vertex weight inequality and fixed vertex constraint are removed by using the idea of equilibrium term and elimination method, and further the integer programming model is relaxed into an easy-to-solve unconstrained optimization problem with the help of trigonometric properties, which is solved by the accelerated subspace minimization conjugate gradient method. Secondly, for the obtained optimal solution of the unconstrained problem the hyperplane rounding algorithm is used for random perturbation to obtain the initial feasible partition of the graph partitioning problem, and further a heuristic algorithm is designed to locally improve the quality of the initial feasible partition solution to finally obtain the approximate optimal solution of the graph partitioning problem. Finally, the numerical performance of our method is tested using a graph partition problem under knapsack constraints and a hypergraph partitioning problem derived from a circuit partitioning industry case study. The experimental results show that the proposed method can effectively solve the multi-constraint graph partitioning problem.

Although the accelerated Subspace Minimization Conjugate Gradient (SMCG) method demonstrates potential for solving hypergraph partitioning problems, its application still faces several challenges, including:

(1) Non-Convexity of the Objective Function: The objective function in hypergraph partitioning is typically non-convex, which may result in multiple local optima. The accelerated SMCG method may converge to a local optimum, posing a challenge for designing global optimization strategies, such as integrating stochastic optimization or multi-start search techniques.

(2) Handling Balanced Constraints: Hypergraph partitioning often requires the results to satisfy balanced constraints (e.g., each subgraph should contain a similar number of nodes). This increases the complexity of the problem. A key challenge lies in effectively incorporating these balanced constraints during subspace minimization while ensuring the algorithm's convergence.

(3) Sparsity and Irregularity of Hypergraph Structures: Hypergraphs are often sparse and irregular, which may lead to instability in subspace construction. Extracting meaningful subspace information from such irregular and sparse hypergraph structures remains a significant challenge.

(4) Parallelization and Distributed Computing: For large-scale hypergraphs, single-machine computation may be insufficient, necessitating distributed computing. Designing efficient parallelization strategies to enable the SMCG method to operate effectively in distributed environments is another critical challenge.

In summary, the Subspace Minimization Conjugate Gradient  method provides an efficient approach for solving hypergraph partitioning problems. However, challenges related to non-convex objective functions, balanced constraints, and other aspects still require further investigation. Future work could focus on integrating adaptive optimization techniques, hybrid strategies, and distributed computing to further enhance the algorithm's performance and applicability.









\begin{acknowledgements}
We would like to thank Prof Hailong You and Dr Fang Zhang for providing the actual industrial design data.
This research was supported by National Natural Science Foundation of China (No. 12261019).
\end{acknowledgements}

\subsection*{\textbf{Declarations}}
\noindent\textbf{Conflict of interest}
The authors declare no competing interests.

\end{document}